\documentclass[a4paper,10pt]{amsart}

\usepackage[english]{babel}
\usepackage[latin1]{inputenc}
\usepackage{amsmath}
\usepackage{amsthm}
\usepackage{amsfonts}
\usepackage{amssymb, latexsym, graphics, showidx}

\numberwithin{equation}{section}
\newtheorem{de}{Definition}[section]
\newtheorem{thm}{Theorem}[section]
\newtheorem{rem}[thm]{Remark}
\newtheorem{cor}[thm]{Corollary}
\newtheorem{prop}[thm]{Proposition}
\newtheorem{lem}[thm]{Lemma}
\newtheorem*{ack}{Acknowledgments}

\renewcommand{\dim}{\noindent\textbf{Proof.} }
\newcommand{\dims}{\noindent\textbf{Proof} }
\newcommand{\finedim}{{\unskip\nobreak\hfil\penalty50
   \hskip2em\hbox{}\nobreak\hfil\mbox{$\Box$ \qquad}
   \parfillskip=0pt \finalhyphendemerits=0\par\medskip}}
\newcommand{\R}{\mathbb{R}}
\newcommand{\N}{\mathbb{N}}
\newcommand{\Om}{\Omega}
\newcommand{\Oms}{\overline{\Omega}}
\newcommand{\lam}{\lambda}
\newcommand{\al}{\alpha}
\newcommand{\lams}{\overline{\lambda}}
\newcommand{\xs}{\overline{x}}
\newcommand{\ys}{\overline{y}}
\newcommand{\zs}{\overline{z}}

\title[The Neumann problem for singular fully nonlinear
operators]{The Neumann problem for singular fully nonlinear
operators}
\author{Stefania Patrizi}
\address{Unversit\`a di Roma ``La Sapienza'', Dipartimento di Matematica,
Piazzale A.~Moro 2, I-00185 Roma, Italy}
\email{patrizi@mat.uniroma1.it}

\begin{document}
\begin{abstract}We consider the Neumann problem in $C^2$ bounded domains for
fully nonlinear second order operators which are elliptic,
homogenous with lower order terms. Inspired by \cite{bnv}, we
define the concept of principal eigenvalue and we characterize it
through the maximum principle. Moreover, Lipschitz regularity,
uniqueness and existence results for solutions of the Neumann
problem are given.
\end{abstract}
%\date{today}
\keywords{Fully nonlinear operators, Neumann problem, Maximum
Principle, principal eigenvalue, viscosity solutions.}

\maketitle

\section{Introduction}
In this paper we study the maximum principle, principal
eigenvalues, regularity and existence for viscosity solutions of
the Neumann boundary value problem
\begin{equation}\label{sisintro}
\begin{cases}
F(x,Du,D^2u)+b(x)\cdot Du|Du|^\al+(c(x)+\lam)|u|^{\al}u= g(x)  & \text{in} \quad\Om \\
 \langle Du, \overrightarrow{n}(x)\rangle= 0 & \text{on} \quad\partial\Om, \\
 \end{cases}
 \end{equation}
where  $\Om$ is a bounded domain of class $C^2$,
$\overrightarrow{n}(x)$ is the exterior normal to the domain $\Om$
at $x$, $\al>-1$, $\lam\in\R$ and $b,c,g$ are continuous functions
on $\Oms$. $F$ is a fully nonlinear operator that may be singular
at the points where the gradient vanishes. It is defined on
$\overline{\Om}\times\R^N \setminus\{0\}\times \emph{S(N)}$, where
$\emph{S(N)}$ denotes  the space of symmetric matrices on $\R^N$
equipped with the usual ordering, and satisfies the following
homogeneity and ellipticity conditions
\renewcommand{\labelenumi}{(F\arabic{enumi})}
\begin{enumerate}
 \item For all $t\in \R^*,\mu \geq 0,$ $(x,p,X)\in\overline{\Om}\times\R^N \setminus\{0\}\times \emph{S(N)}$
$$F(x,tp,\mu X)=|t|^\al\mu F(x,p,X).$$
 \item There exist $a,A>0$ such that for $x\in\overline{\Om},\, p\in \R^N\setminus\{0\}, M, N\in \emph{S(N)},\,N\geq 0$
 $$ a|p|^\al \text{tr}N\leq F(x,p,M+N)-F(x,p,M)\leq A|p|^\al \text{tr}N.$$
\end{enumerate}
In addition, we will assume on $F$ some H\"{o}lder's continuity
hypothesis that will be made precise in the next section.

In this class of operators one can consider for example
$$F(Du,D^2u)=|Du|^\al \mathcal{M}_{a,A}^{+\atop -}D^2u,$$ $\al>-1$, where $\mathcal{M}_{a,A}^{+\atop -}D^2u$ are the Pucci's operators
(see e.g. \cite{cc}), the p-Laplacian
$$\Delta_pu=\text{div}(|Du|^{p-2}Du),$$ with $\al=p-2$, and
non-variational extensions of the p-Laplacian, depending
explicitly on $x$, like the operator
\begin{equation*}F(x,Du,D^2u)=
|Du|^{q-2}\text{tr}(B_1(x)D^2u)+c_0|Du|^{q-4}\langle
D^2uB_2(x)Du,B_2(x)Du\rangle,\end{equation*} with $\al=q-2$, where
$q>1$, $B_1$ and $B_2$ are $\theta$-H\"{o}lderian functions with
$\theta>\frac{1}{2}$, which send $\overline{\Omega}$ into
\emph{S(N)}, $aI\leq B_1\leq AI$, $-\sqrt{a}I\leq B_2\leq
\sqrt{a}I$ and $c_0>-1$.

The concept of first eigenvalue has been extended to nonlinear
operators which are variational, such as the p-Laplacian with
Dirichlet or Neumann boundary conditions, through the method of
minimization of the Rayleigh quotient, see e.g. \cite{an} and
\cite{li2}. That method uses heavily the variational structure and
cannot be applied to operators which have not this property. An
important step in the study of the eigenvalue problem for
nonlinear operators in non-divergence form was made  by Lions in
\cite{l}. In that paper, using probabilistic and analytic methods,
he showed the existence of principal eigenvalues for the uniformly
elliptic Bellman operator and obtained results for the related
Dirichlet problems. Very recently, many authors, inspired by the
famous work of Berestycki, Nirenberg and Varadhan \cite{bnv}, have
de\-ve\-lo\-ped an eigenvalue theory for fully nonlinear operators
which are non-variational. Issues similar to those of this paper
have been studied for the Dirichlet problem by Birindelli and
Demengel in \cite{bd}. They assume slightly less general structure
conditions on $F$, but on the other hand, some of their results
can be applied to degenerate elliptic equations. The case $\al=0$
has been treated by Quaas \cite{q} and Busca, Esteban and Quaas
\cite{beq} for the Pucci's operators. Their results have been
extended to more general fully nonlinear convex uniformly elliptic
operators in \cite{qs} by Quaas and Sirakov. See also the work of
Ishii and Yoshimura \cite{iy} for non-convex operators. All these
articles treat Dirichlet boundary conditions.

The techniques of this paper, although partly taken by the
previous mentioned articles, use ad hoc test functions depending
on the distance function from the boundary of the domain which are
suitable for the Neumann boundary conditions.

Comparison principles and the existence results for the Neumann
problem have been investigated by Ishii in \cite{i} and Barles in
\cite{b2} and \cite{b} for degenerate elliptic operators
$\mathcal{G}(x,u,Du,D^2u)$ modeled on the Isaacs ones or on the
stationary operator associated to the Mean Curvature Equation. In
all these papers a fundamental  assumption is the monotonicity of
$\mathcal{G}(x,r,p,X)$ with respect to $r$. For the p-Laplace with
the zero order term $c(x)|u|^{p-2}u$, $c\leq 0$ and $c\not\equiv
0$ and the pure Neumann boundary condition, the comparison
principle can be showed through variational techniques, like in
the Dirichlet case, see e.g. \cite{li}.

We denote
\begin{equation*}G(x,u,Du,D^2u):=
F(x,Du,D^2u)+b(x)\cdot Du|Du|^\al+c(x)|u|^{\al}u.\end{equation*}It
is important to remark that $G$ is homogenous  and
non-variational. Following the ideas of \cite{bnv}, we define the
principal eigenvalue as
\begin{equation}\label{primoaut}\begin{split}\lams:=\sup\{&\lam\in
\R\;|\;\exists\,v>0 \text{ bounded viscosity solution of }\\&
G(x,v,Dv,D^2v)+\lam v^{\al+1}\leq 0 \text{ in } \Om,\langle
Dv,\overrightarrow{n}\rangle \geq 0\text{ on }\partial \Om \}.
\end{split}\end{equation}
$\lams$ is well defined since the above set is not empty; indeed,
$-|c|_\infty$ belongs to it, being $v(x)\equiv1$ a corresponding
supersolution. Furthermore it is an interval because if $\lam$
belongs to it then so does any $\lam'<\lam$.

One of the scope of this work is to prove that $\lams$ is an
"eigenvalue" for $-G$ which admits a positive "eigenfunction", in
the sense that there exists $\phi>0$ solution of
\begin{equation*}
\begin{cases}
G(x,\phi,D\phi,D^2\phi)+\lams \phi^{\al+1}= 0  & \text{in} \quad\Om \\
 \langle D\phi, \overrightarrow{n}(x)\rangle= 0 & \text{on} \quad\partial\Om. \\
 \end{cases}
 \end{equation*}
Moreover, $\lams$ can be characterized as the supremum of those
$\lam$ for which the ope\-ra\-tor $G(x,u,Du,D^2u)+\lam |u|^\al u$
with the Neumann boundary condition satisfies the maximum
principle. In particular $\lams$ is the least "eigenvalue" to
which there cor\-re\-spond "eigenfunctions" positive somewhere.
These results are applied to obtain existence and uniqueness for
the boundary value problem \eqref{sisintro}.

The paper is organized as follows. In the next section we give
assumptions and define the concept of solution. In Section 3 we
establish a Lipschitz regularity result for viscosity solutions of
\eqref{sisintro}. The Section 4 is devoted to the study of the
maximum principle for subsolutions of \eqref{sisintro}.  In
Section 4.1 we show that it holds (even for more general boundary
conditions) for $G(x,u,Du,D^2u)$ if $c(x)\leq 0$ and $c\not\equiv
0$, see Theorem \ref{macpcless0}. One of the main result of the
paper is that the maximum principle holds for $G(x,u,Du,D^2u)+\lam
|u|^\al u$ for any $\lam<\lams$, as we show in Theorem
\ref{maxpneum} of Section 4.2. In particular it holds for
$G(x,u,Du,D^2u)$ if $\lams>0$. It is natural to wonder if the
result of Theorem \ref{maxpneum} is stronger than that of Theorem
\ref{macpcless0}; indeed if $c\equiv 0$, one has $\lams=0$. A
positive answer is given in Section 4.3, where we construct an
explicit example of a bounded positive viscosity solution of
$G(x,v,Dv,D^2v)+\lam v^{\al+1}\leq 0$ in $\Om$, $\langle
Dv,\overrightarrow{n}\rangle \geq 0$ on $\partial \Om$, $\lam>0$,
with $c(x)$ changing sign. The existence of such $v$ implies, by
definition, $\lams>0$. Finally, in Section 5 we show some
existence and comparison theorems.

For fully nonlinear operators it is possible to define another
principal eigenvalue
\begin{equation*}\begin{split}\underline{\lam}:=\sup\{&\lam\in
\R\;|\;\exists\,u<0 \text{ bounded viscosity solution of }\\&
G(x,u,Du,D^2u)+\lam |u|^{\al}u\geq 0 \text{ in } \Om,\langle
Du,\overrightarrow{n}\rangle \leq 0\text{ on }\partial \Om \}.
\end{split}\end{equation*}
If $F(x,p,X)=-F(x,p,-X)$ then $\lams=\underline{\lam}$, otherwise
$\lams$ may be different from $\underline{\lam}$.

The classical assumption which guarantees the solvability of the
Neumann pro\-blem \eqref{sisintro} is $c<0$ in $\Oms$. We show
that the right hypothesis for any right-hand side is the
positivity of the two principal eigenvalues.

\section{Assumptions and definitions}
We assume that the operator $F:\overline{\Om}\times\R^N
\setminus\{0\}\times \emph{S(N)}\rightarrow\R$ satisfies the
hy\-po\-the\-sis (F1) and (F2) given in the introduction and the
following H\"{o}lder's continuity conditions
\begin{itemize}
\item [(F3)] There exist $C_1>0$ and $\theta\in (\frac{1}{2},1]$ such
that for all $x,y\in
\Oms,\,p\in\R^N\setminus\{0\},\,X\in\emph{S(N)}$
$$|F(x,p,X)-F(y,p,X)|\leq
C_1|x-y|^\theta|p|^\al\|X\|.$$
\item [(F4)] There exist $C_2>0$ and $\nu\in
(\frac{1}{2},1]$ such that for all $x\in\Oms$,
$p\in\R^N\setminus\{0\}$, $p_0\in\R^N$, $|p_0|\leq \frac{|p|}{2}$,
$X\in\emph{S(N)}$
$$|F(x,p+p_0,X)-F(x,p,X)|\leq C_2|p|^{\al-\nu}|p_0|^\nu\|X\|.$$
\end{itemize}
Here and in what follows we fix the norm $\|X\|$ in $\emph{S(N)}$
by setting
$\|X\|=\sup\{|X\xi|\,|\,\xi\in\R^N,\,|\xi|\leq1\}=\sup\{|\lam|\,:\,
\lam \text{ is an eigenvalue of }X\}.$

The domain $\Om$ is supposed to be bounded and of class $C^2$. In
particular, it satisfies the interior sphere condition and the
uniform exterior sphere condition, i.e.,
\begin{enumerate}
\item [($\Om1$)]For
 each $x\in\partial \Om$ there exist $R>0$ and $y\in \Om$ for
 which $|x-y|=R$ and $B(y,R)\subset \Om$.
\item[($\Om2$)] There exists $r>0$ such that  $B(x+r\overrightarrow{n}(x),r)\cap \Om =\emptyset$  for any $x\in \partial\Om.$\end{enumerate}
From the property ($\Om2$) it follows that
\begin{equation}\label{sferaest}
\langle
y-x,\overrightarrow{n}(x)\rangle\leq\frac{1}{2r}|y-x|^2\quad\text{for
}x\in\partial\Om \text{ and } y\in\Oms.\end{equation} Moreover,
the $C^2$-regularity of $\Om$ implies the existence of a
neighborhood of $\partial\Om$ in $\Oms$ on which the distance from
the boundary
 $$d(x):=\inf\{|x-y|, y\in\partial\Om\},\quad x\in\Oms$$ is of class $C^2$. We still denote by $d$ a
$C^2$ extension of the distance function to the whole $\Oms$.
Without loss of generality we can assume that $|Dd(x)|\leq1$ in
$\Oms$.

As in \cite{bd}, here we adopt the notion of viscosity solution
for \eqref{sisintro} adapted to our context. We denote by
$USC(\overline{\Om})$ the set of upper semicontinuous functions on
$\overline{\Om}$ and by $LSC(\overline{\Om})$ the set of lower
semicontinuous functions on $\overline{\Om}$. Let
$g:\Oms\rightarrow \R$ and $B:\partial\Om\times\R\times
\R^N\rightarrow\R$.
\begin{de}Any function $u\in  USC(\overline{\Om})$  (resp., $u\in  LSC(\overline{\Om})$) is called \emph{viscosity subsolution}
(resp., \emph{supersolution}) of
\begin{equation*}
\begin{cases}
 G(x,u,Du,D^2u)= g(x)  & \text{in} \quad\Om \\
 B(x,u,Du)= 0 & \text{on} \quad\partial\Om, \\
 \end{cases}
 \end{equation*}
 if the following conditions hold
\begin{itemize}
\item[(i)] For every $x_0\in \Om$, for all $\varphi\in C^2(\overline{\Om})$, such that $u-\varphi$ has a local maximum (resp., minimum)
on $x_0$ and $D\varphi(x_0)\neq0$, one has
$$G(x_0,u(x_0),D\varphi(x_0),D^2\varphi(x_0))\geq \,(\text{resp., } \leq\,)\,\,g(x_0).$$ If
 $u\equiv k=$const. in a neighborhood of $x_0$, then
$$c(x_0)|k|^\al k\geq\, (\text{resp., } \leq\,)\,\,g(x_0).$$
\item[(ii)]For every $x_0\in \partial\Om$, for all $\varphi\in C^2(\overline{\Om})$, such that $u-\varphi$ has a local maximum (resp., minimum)
on $x_0$ and $D\varphi(x_0)\neq0$, one has
$$(-G(x_0,u(x_0),D\varphi(x_0),D^2\varphi(x_0))+g(x_0))\wedge
B(x_0,u(x_0),D\varphi(x_0))\leq 0$$
(resp.,$$(-G(x_0,u(x_0),D\varphi(x_0),D^2\varphi(x_0))+g(x_0))\vee
B(x_0,u(x_0),D\varphi(x_0))\geq 0).$$ If $u\equiv k=$const. in a
neighborhood of $x_0$ in $\overline{\Om}$, then
$$(-c(x_0)|k|^\al k+g(x_0))\wedge B(x_0,k,0)\leq 0$$\quad (resp., $$(-c(x_0)|k|^\al k+g(x_0))\vee B(x_0,k,0)\geq 0).$$
\end{itemize}
A \emph{viscosity solution} is a continuous function which is both
a subsolution and a supersolution.
\end{de}
For a detailed presentation of the theory of viscosity solutions
and of the boun\-da\-ry conditions in the viscosity sense, we
refer the reader to \cite{cil}.

We call strong viscosity subsolutions (resp., supersolutions) the
viscosity subsolutions (resp., supersolutions) that satisfy
$B(x,u,Du)\leq $ (resp., $\geq$) 0 in the viscosity sense for all
$x\in\partial\Om$. If $\lam\rightarrow
B(x,r,p-\lam\overrightarrow{n})$ is non-increasing in $\lam\geq0$,
then classical subsolutions (resp., supersolutions) are strong
viscosity subsolutions (resp., supersolutions), see \cite{cil}
Proposition 7.2.

In the above definition the test functions can be substituted by
the elements of the set $\overline{J}^{2,+}u(x_0)$ when $u$ is a
subsolution and $ \overline{J}^{2,-}u(x_0)$ when $u$ is a
supersolution. For non-singular operators the definitions reduce
to the standard ones, see \cite{cil}.

\section{Lipschitz continuity of viscosity solutions}
\begin{thm}\label{regolaritathm}Let $\Om$ be a bounded domain of class $C^2$. Suppose that $F$ satisfies (F2)-(F4) and that $b$, $c$, $g$
are bounded in $\Om$. If $u\in C(\Oms)$ is a viscosity solution of
\begin{equation*}
\begin{cases}
 F(x,Du,D^2u)+b(x)\cdot Du|Du|^\al+c(x)|u|^\al u=  g(x) & \text{in} \quad\Om \\
 \langle Du, \overrightarrow{n}(x)\rangle= 0 & \text{on} \quad\partial\Om, \\
 \end{cases}
 \end{equation*}then
$$ |u(x)-u(y)|\leq C_0|x-y|\quad \forall x,y\in\Oms,$$
where $C_0$ depends on
 $\Om,\,N,\,a,\,A,\,\theta,\,\nu,\,C_1,\,C_2,\,|b|_\infty,\,|c|_\infty,\,|g|_\infty,$ and
 $|u|_\infty$.
 \end{thm}
 The Theorem is an immediate consequence of the next lemma.
  To prove the lemma we adopt the technique used
in Proposition III.1 of \cite{il} for Dirichlet pro\-blems, that
we modify taking test functions which depend on $d(x)$.

The lemma plays a key role also in the proof of Theorem
\ref{maxpneum} in the next section.

\begin{lem}\label{regolarita}Assume the hypothesis of Theorem \ref{regolaritathm} and suppose that $g$ and $h$ are bounded functions.
Let $u\in USC(\Oms)$ be a viscosity subsolution of
\begin{equation*}
\begin{cases}
 F(x,Du,D^2u)+b(x)\cdot Du|Du|^\al+c(x)|u|^\al u=  g(x) & \text{in} \quad\Om \\
 \langle Du, \overrightarrow{n}(x)\rangle= 0 & \text{on} \quad\partial\Om, \\
 \end{cases}
 \end{equation*}
and $v\in LSC(\Oms)$ a viscosity supersolution of
\begin{equation*}
\begin{cases}
 F(x,Dv,D^2v)+b(x)\cdot Dv|Dv|^\al+c(x)|v|^\al v=  h(x) & \text{in} \quad\Om \\
 \langle Dv, \overrightarrow{n}(x)\rangle= 0 & \text{on} \quad\partial\Om, \\
 \end{cases}
 \end{equation*}
with $u$ and $v$  bounded, or $v\geq 0$ and bounded. If
$m=\max_{\Oms}(u-v)\geq0$,
 then there exists $C_0>0$ such that
\begin{equation}\label{stimau-v1}u(x)-v(y)\leq m +C_0|x-y|\quad \forall
x,y\in\Oms,\end{equation}
 where $C_0$ depends on
 $\Om,\,N,\,a,\,A,\,\theta,\,\nu,\,C_1,\,C_2,\,|b|_\infty,\,|c|_\infty,\,|g|_\infty,\,|h|_\infty,\,|v|_\infty,\,m$ and
 $|u|_\infty$ or $\sup_{\Oms}u$.
 \end{lem}
 \dim
We set
$$\Phi(x)=MK|x|-M(K|x|)^2,$$and
$$\varphi (x,y)=m+e^{-L(d(x)+d(y))}\Phi(x-y),$$ where
$L$ is a fixed number greater than $2/(3r)$ with $r$ the radius in
the condition ($\Om2$) and $K$ and $M$ are two positive constants
to be chosen later. If $K|x|\leq \frac{1}{4}$, then
\begin{equation}\begin{split}\label{phimagg}\Phi(x)&=MK|x|-M(K|x|)^2\geq
\frac{3}{4}MK|x|.\end{split}\end{equation} We define
$$\Delta _K:=\left\{(x,y)\in \R^N\times\R^N|\, |x-y|\leq\frac{1}{4K}\right\}.$$ We
fix $M$ such that
\begin{equation}\label{M}
\max_{\Oms^{\,2}}(u(x)-v(y))\leq
m+e^{-2Ld_0}\frac{M}{8},\end{equation} where
$d_0=\max_{x\in\Oms}d(x),$ and we claim that taking $K$ large
enough, one has
$$u(x)-v(y)-\varphi(x,y)\leq 0\quad\text{for }(x,y)\in
\Delta_K\cap\Oms^2.$$ In this case \eqref{stimau-v1} is proven. To
show the last inequality we suppose by contradiction that for some
$(\xs,\ys)\in\Delta_K\cap \Oms^2$
$$u(\xs)-v(\ys)-\varphi(\xs,\ys)=\max_{\Delta_K\cap
\Oms\,^2}(u(x)-v(y)-\varphi(x,y))>0.$$ Here we have dropped the
dependence of $\xs,\,\ys$ on $K$ for simplicity of notations.

Observe that if $v\geq0$, since from \eqref{phimagg} $\Phi(x-y)$
is non-negative in $\Delta_K$ and $m\geq 0$, one has $u(\xs)>0$.

Clearly $\xs\neq \ys$. Moreover the point $(\xs,\ys)$ belongs to
$\text{int}(\Delta_K)\cap\Oms^2$. Indeed, if $|x-y|=\frac{1}{4K}$,
by \eqref{M} and \eqref{phimagg} we have
\begin{equation*}
u(x)-v(y)\leq m +e^{-2Ld_0} \frac{M}{8}\leq
m+e^{-L(d(x)+d(y))}\frac{1}{2}MK|x-y|\leq\varphi(x,y).\end{equation*}
Since $\xs\neq \ys$ we can compute the derivatives of $\varphi$ in
$(\xs,\ys)$ obtaining
\begin{equation*}\begin{split}D_x\varphi(\xs,\ys)&=e^{-L(d(\xs)+d(\ys))}MK\Big\{-L|\xs-\ys|(1-K|\xs-\ys|)Dd(\xs)\\&+
(1-2K|\xs-\ys|)\frac{(\xs-\ys)}{|\xs-\ys|}\Big\},\end{split}\end{equation*}
\begin{equation*}\begin{split}D_y\varphi(\xs,\ys)&=e^{-L(d(\xs)+d(\ys))}MK\Big\{-L|\xs-\ys|(1-K|\xs-\ys|)Dd(\ys)\\&-
(1-2K|\xs-\ys|)\frac{(\xs-\ys)}{|\xs-\ys|}\Big\}.\end{split}\end{equation*}
Observe that for large $K$
\begin{equation}\label{dxphistima}\begin{split}e^{-2Ld_0}\frac{MK}{4}&\leq e^{-L(d(\xs)+d(\ys))}MK\left(\frac{1}{2}-L|\xs-\ys|\right)\leq
|D_x\varphi(\xs,\ys)|,|D_y\varphi(\xs,\ys)|\\&\leq
2MK.\end{split}\end{equation} Using \eqref{sferaest}, if
$\xs\in\partial \Om$ we have
\begin{equation*}\begin{split}
&\langle D_x\varphi(\xs,\ys),\overrightarrow{n}(\xs)\rangle\\& =
e^{-Ld(\ys)}MK\Big\{L|\xs-\ys|(1-K|\xs-\ys|)+(1-2K|\xs-\ys|)\langle\frac{(\xs-\ys)}{|\xs-\ys|},\overrightarrow{n}(\xs)\rangle\Big\}\\&\geq
e^{-Ld(\ys)}MK\Big\{\frac{3}{4}L|\xs-\ys|-(1-2K|\xs-\ys|)\frac{|\xs-\ys|}{2r}\Big\}\\&\geq
\frac{1}{2}e^{-Ld(\ys)}MK|\xs-\ys|\left(\frac{3}{2}L-\frac{1}{r}\right)>0,
\end{split}\end{equation*}since
$\xs\neq\ys$ and $L>2/(3r)$. Similarly,
 if $\ys\in\partial\Om$
\begin{equation*}
\langle-D_y\varphi(\xs,\ys),\overrightarrow{n}(\ys)\rangle\leq
\frac{1}{2}e^{-Ld(\xs)}MK|\xs-\ys|\left(-\frac{3}{2}L+\frac{1}{r}\right)<0.\end{equation*}
In view of definition of sub and supersolution, we conclude that
$$G(\xs,u(\xs),D_x\varphi(\xs,\ys),X)\geq g(\xs)\quad\text{if }(D_x\varphi(\xs,\ys),X)\in
\overline{J}^{2,+}u(\xs),$$
$$G(\ys,v(\ys),-D_y\varphi(\xs,\ys),Y)\leq h(\ys)\quad\text{if }(-D_y\varphi(\xs,\ys),Y)\in \overline{J}^{2,-}v(\ys).$$ Since
$(\xs,\ys)\in\text{int}\Delta_K\cap\Oms\,^2$, it is a local
maximum point of $u(x)-v(y)-\varphi(x,y)$ in $\Oms\,^2$. Then
applying Theorem 3.2 in \cite{cil}, for every $\epsilon>0$ there
exist $X,Y\in \emph{S(N)}$ such that $ (D_x\varphi(\xs,\ys),X)\in
\overline{J}\,^{2,+}u(\xs),\,(-D_y\varphi(\xs,\ys),Y)\in
\overline{J}\,^{2,-}v(\ys)$ and
\begin{equation}\label{tm2ishii}\left(%
\begin{array}{cc}
  X & 0 \\
  0 & -Y \\
\end{array}%
\right)\leq D^2(\varphi(\xs,\ys))+\epsilon
(D^2(\varphi(\xs,\ys)))^2.
\end{equation}
Now we want to estimate the matrix on the right-hand side of the
last inequality.
\begin{equation*}\begin{split}D^2\varphi(\xs,\ys)&=\Phi(\xs-\ys)D^2(e^{-L(d(\xs)+d(\ys))})+D(e^{-L(d(\xs)+d(\ys))})\otimes
D(\Phi(\xs-\ys))\\&+D(\Phi(\xs-\ys))\otimes
D(e^{-L(d(\xs)+d(\ys))})+e^{-L(d(\xs)+d(\ys))}D^2(\Phi(\xs-\ys)).\end{split}\end{equation*}We
set $$A_1:=\Phi(\xs-\ys)D^2(e^{-L(d(\xs)+d(\ys))}),$$
$$A_2:=D(e^{-L(d(\xs)+d(\ys))})\otimes
D(\Phi(\xs-\ys))+D(\Phi(\xs-\ys))\otimes
D(e^{-L(d(\xs)+d(\ys))}),$$
$$A_3:=e^{-L(d(\xs)+d(\ys))}D^2(\Phi(\xs-\ys)).$$Observe that
\begin{equation}\label{a1}A_1\leq CK|\xs-\ys|\left(%
\begin{array}{cc}
  I & 0 \\
  0 & I \\
\end{array}%
\right).\end{equation} Here and henceforth C denotes various
positive constants independent of $K$.

For $A_2$ we have the following estimate
\begin{equation}\label{a2}A_2\leq
CK\left(%
\begin{array}{cc}
  I & 0 \\
  0 & I \\
\end{array}%
\right)
+CK\left(%
\begin{array}{cc}
  I & -I \\
  -I & I \\
\end{array}%
\right).\end{equation} Indeed for $\xi,\,\eta\in\R^N$ we compute
\begin{equation*}\begin{split}\langle
A_2(\xi,\eta),(\xi,\eta)\rangle&=2Le^{-L(d(\xs)+d(\ys))}\{\langle
Dd(\xs)\otimes D\Phi(\xs-\ys)(\eta-\xi),\xi\rangle\\&+\langle
Dd(\ys)\otimes D\Phi(\xs-\ys)(\eta-\xi),\eta\rangle\}\leq
CK(|\xi|+|\eta|)|\eta-\xi|\\&\leq
CK(|\xi|^2+|\eta|^2)+CK|\eta-\xi|^2.\end{split}\end{equation*}

Now we consider $A_3$. The matrix $D^2(\Phi(\xs-\ys))$ has the
form
$$D^2(\Phi(\xs-\ys))=\left(%
\begin{array}{cc}
  D^2\Phi(\xs-\ys) & - D^2\Phi(\xs-\ys) \\
  - D^2\Phi(\xs-\ys) &  D^2\Phi(\xs-\ys) \\
\end{array}%
\right),$$and the Hessian matrix of $\Phi(x)$ is
\begin{equation}\label{hessianphi}D^2\Phi(x)=\frac{MK}{|x|}\left(I-\frac{x\otimes
x}{|x|^2}\right)-2MK^2I.\end{equation} If we choose
\begin{equation*}\epsilon=\frac{|\xs-\ys|}{2MKe^{-L(d(\xs)+d(\ys))}},\end{equation*}
then we have the following estimates
$$\epsilon A_1^2\leq
CK|\xs-\ys|^3I_{2N},\quad \epsilon A_2^2\leq CK|\xs-\ys|I_{2N},$$
\begin{equation}\label{aprodotti}\begin{split} \epsilon (A_1A_2+A_2A_1)\leq
CK|\xs-\ys|^2I_{2N},\end{split}\end{equation}
$$\epsilon (A_1A_3+A_3A_1)\leq CK|\xs-\ys|I_{2N},\quad \epsilon
(A_2A_3+A_3A_2)\leq CKI_{2N},$$
 where $I_{2N}:=\left(%
\begin{array}{cc}
  I & 0 \\
  0 & I \\
\end{array}%
\right)$. Then using \eqref{a1}, \eqref{a2},
\eqref{aprodotti} and observing that $$(D^2(\Phi(\xs-\ys)))^2=\left(%
\begin{array}{cc}
  2(D^2\Phi(\xs-\ys))^2 & - 2(D^2\Phi(\xs-\ys))^2 \\
  - 2(D^2\Phi(\xs-\ys))^2 &  2(D^2\Phi(\xs-\ys))^2 \\
\end{array}%
\right),$$from \eqref{tm2ishii} we conclude that
$$\left(%
\begin{array}{cc}
  X & 0 \\
  0 & -Y \\
\end{array}%
\right)\leq O(K)\left(%
\begin{array}{cc}
  I & 0 \\
  0 & I \\
\end{array}%
\right)+\left(%
\begin{array}{cc}
  B & -B \\
  -B & B \\
\end{array}%
\right),$$ where
\begin{equation*}B=CKI+e^{-L(d(\xs)+d(\ys))}\left[D^2\Phi(\xs-\ys)+\frac{|\xs-\ys|}{MK}
(D^2\Phi(\xs-\ys))^2\right].\end{equation*}The last inequality can
be rewritten as follows
$$\left(%
\begin{array}{cc}
  \widetilde{X} & 0 \\
  0 & -\widetilde{Y} \\
\end{array}%
\right)\leq\left(%
\begin{array}{cc}
  B & -B \\
  -B & B \\
\end{array}%
\right),$$ with $\widetilde{X}=X-O(K)I$ and
$\widetilde{Y}=Y+O(K)I.$

Now we want to get a good estimate for
tr($\widetilde{X}-\widetilde{Y}$), as in \cite{il}. For that aim
let
$$0\leq P:=\frac{(\xs-\ys)\otimes (\xs-\ys)}{|\xs-\ys|^2}\leq I.$$
Since $\widetilde{X}-\widetilde{Y}\leq 0$ and
$\widetilde{X}-\widetilde{Y}\leq 4B,$ we have
$$\text{tr}(\widetilde{X}-\widetilde{Y})\leq \text{tr}(P(\widetilde{X}-\widetilde{Y}))\leq 4 \text{tr}(PB).$$ We
have to compute tr($PB$). From \eqref{hessianphi}, observing that
the matrix $(1/|x|^2)x\otimes x$ is idempotent, i.e.,
$[(1/|x|^2)x\otimes x]^2=(1/|x|^2)x\otimes x$, we compute
$$(D^2\Phi(x))^2=\frac{M^2K^2}{|x|^2}(1-4K|x|)\left(I-\frac{x\otimes
x}{|x|^2}\right)+4M^2K^4I.$$ Then, since $\text{tr}P=1$ and
$4K|\xs-\ys|\leq1$, we have
\begin{equation*}\begin{split}\text{tr}(PB)&=CK+e^{-L(d(\xs)+d(\ys))}(-2MK^2+4MK^3|\xs-\ys|)
\\&\leq CK-e^{-L(d(\xs)+d(\ys))}MK^2<0, \end{split}\end{equation*}for
large $K$. This gives
$$|\text{tr}(\widetilde{X}-\widetilde{Y})|=-\text{tr}(\widetilde{X}-\widetilde{Y})\geq
4e^{-L(d(\xs)+d(\ys))}MK^2-4CK\geq CK^2,$$
 for large $K$. Since $\|B\|\leq \frac{CK}{|\xs-\ys|},$ we have
\begin{equation*}\begin{split}\|B\|^{\frac{1}{2}}|\text{tr}(\widetilde{X}-\widetilde{Y})|^{\frac{1}{2}}&\leq
\left(\frac{CK}{|\xs-\ys|}\right)^{\frac{1}{2}}|\text{tr}(\widetilde{X}-\widetilde{Y})|^{\frac{1}{2}}
\leq
\frac{C}{K^{\frac{1}{2}}|\xs-\ys|^{\frac{1}{2}}}|\text{tr}(\widetilde{X}-\widetilde{Y})|.
\end{split}\end{equation*}The Lemma III.I in \cite{il} ensures the
existence of a universal constant $C$ depending only on $N$ such
that $$\|\widetilde{X}\|, \|\widetilde{Y}\|\leq
C\{|\text{tr}(\widetilde{X}-\widetilde{Y})|+\|B\|^{\frac{1}{2}}|\text{tr}(\widetilde{X}-\widetilde{Y})|^{\frac{1}{2}}\}.$$
Thanks to the above estimates we can conclude that
\begin{equation}\label{normaxy}\|\widetilde{X}\|,\,\|\widetilde{Y}\|\leq
C|\text{tr}(\widetilde{X}-\widetilde{Y})|\left(1+\frac{1}{K^{\frac{1}{2}}|\xs-\ys|^{\frac{1}{2}}}\right).\end{equation}

Now, using the assumptions (F2), (F3) and (F4) concerning $F$, the
definition of $\widetilde{X}$ and $\widetilde{Y}$ and the fact
that $u$ and $v$ are respectively sub and supersolution we
compute\begin{equation*}\begin{split}g(\xs)-c(\xs)|u(\xs)|^\al
u(\xs)&\leq F(\xs,D_x\varphi,X)+b(\xs)\cdot
D_x\varphi|D_x\varphi|^\al\\&\leq F(\xs,D_x\varphi,\widetilde{X})
+|D_x\varphi|^\al O(K)+b(\xs)\cdot D_x\varphi|D_x\varphi|^\al\\&
\leq
F(\ys,-D_y\varphi,\widetilde{Y})+C_1|\xs-\ys|^\theta|D_x\varphi|^\al\|\widetilde{X}\|\\&
+CK^\nu|\xs-\ys|^{\nu}|D_x\varphi|^{\al-\nu}
\|\widetilde{X}\|+a|D_y\varphi|^\al\text{tr}(\widetilde{X}-\widetilde{Y})\\&+|D_x\varphi|^\al
O(K)+b(\xs)\cdot D_x\varphi|D_x\varphi|^\al\\&\leq b(\ys)\cdot
D_y\varphi|D_y\varphi|^\al -c(\ys)|v(\ys)|^\al v(\ys)+h(\ys)
\\&+C_1|\xs-\ys|^\theta|D_x\varphi|^\al\|\widetilde{X}\|+CK^\nu|\xs-\ys|^{\nu}|D_x\varphi|^{\al-\nu}
\|\widetilde{X}\|\\&+a|D_y\varphi|^\al\text{tr}(\widetilde{X}-\widetilde{Y})+|D_y\varphi|^\al\vee|D_x\varphi|^\al
O(K)\\&+b(\xs)\cdot
D_x\varphi|D_x\varphi|^\al.\end{split}\end{equation*} From this
inequalities, using \eqref{dxphistima}, \eqref{normaxy} and the
fact that $\theta,\nu>\frac{1}{2}$ we get
\begin{equation*}\begin{split}
&g(\xs)-h(\ys)-c(\xs)|u(\xs)|^\al u(\xs)+c(\ys)|v(\ys)|^\al
v(\ys)\leq |D_y\varphi|^\al\vee|D_x\varphi|^\al[
a\text{tr}(\widetilde{X}-\widetilde{Y})\\&+C_1|\xs-\ys|^\theta\|\widetilde{X}\|+C|\xs-\ys|^{\nu}\|\widetilde{X}\|+O(K)]
\leq C
K^\al[a\text{tr}(\widetilde{X}-\widetilde{Y})+o(|\text{tr}(\widetilde{X}-\widetilde{Y})|)].
\end{split}\end{equation*}
If both $u$ and $v$ are bounded, then the first member in the last
inequalities is bounded from below by
$-|g|_\infty-|h|_\infty-|c|_\infty(|u|_\infty^{\al+1}+|v|_\infty^{\al+1})$.
Otherwise, if $v$ is non-negative and bounded, then $u(\xs)\geq 0$
and that quantity is greater than
$-|g|_\infty-|h|_\infty-|c|_\infty (\sup u)^{\al+1}-|c|_\infty
|v|_\infty^{\al+1}$. On the other hand, the last member goes to
$-\infty$ as $K\rightarrow+\infty$, hence taking $K$ large enough
we obtain a contradiction and this concludes the proof.\finedim
\begin{rem}\label{stimadirich}{\em If $F$ satisfies (F2) and (F3),  $u$ is a subsolution of $G(x,u,Du,D^2u)=g$,
$v$ is a supersolution of $G(x,v,Dv,D^2v)=h$ in $\Om$, $u\leq v$
on $\partial\Om$ and $m>0$ then the estimate \eqref{stimau-v1}
still holds for any $x,y\in\Om$. To prove this define
$\varphi=m+MK|x|-M(K|x|)^2$ and follow the proof of Lemma
\ref{regolarita}.}
\end{rem}

Since the Lipschitz constant of the solution depends only on its
bound, on the bound of $g$ and on the structural constants, an
immediate consequence of Theorem \ref{regolaritathm} is the
following compactness criterion that will be useful in the last
section.
\begin{cor}\label{corcomp}Assume the hypothesis of Theorem \ref{regolaritathm} on $\Om$, $F$ and $b$. Suppose that $(g_n)_{n}$ is a
sequence of continuous and uniformly bounded functions and
$(u_n)_n$ is a sequence of uniformly bounded viscosity solutions
of
\begin{equation*}
\begin{cases}
 F(x,Du_n,D^2u_n)+b(x)\cdot Du_n|Du_n|^\al=  g_n(x) & \text{in} \quad\Om \\
 \langle Du_n, \overrightarrow{n}(x)\rangle= 0 & \text{on} \quad\partial\Om. \\
 \end{cases}
 \end{equation*}Then the sequence $(u_n)_n$ is relatively compact in
 $C(\Oms)$.
 \end{cor}
 \section{The Maximum Principle and the principal eigenvalues}
\noindent We say that the operator $ G(x,u,Du,D^2u)$ with the
Neumann boundary condition satisfies the maximum principle if
whenever $u\in USC(\Oms)$ is a viscosity subsolution of
\begin{equation*}
\begin{cases}
 G(x,u,Du,D^2u)= 0  & \text{in} \quad\Om \\
 \langle Du,\overrightarrow{n}(x)\rangle= 0 & \text{on} \quad\partial\Om, \\
 \end{cases}
 \end{equation*}then $u\leq 0$ in $\Oms$.

We first prove that the maximum principle holds under the
classical assumption $c\leq 0$, also for domain which are not of
class $C^2$ and with more general boundary conditions. Then we
show that the operator $G(x,u,Du,D^2u)+\lam |u|^{\al}u$ with the
Neumann boundary condition satisfies the maximum principle for any
$\lam<\lams$. This is the best result that one can expect, indeed,
as we will see in the last section, $\lams$ admits a positive
eigenfunction which provides a counterexample to the maximum
principle for $\lam\geq\lams$.

Finally, we give an example of $c(x)$ which changes sign in $\Om$
and such that the associated principal eigenvalue $\lams$ is
positive.

 \subsection{The case $c(x)\leq 0$}
In this subsection we assume that $\Om$ is of class $C^1$ and
satisfies the interior sphere condition ($\Om$1). We need the
comparison principle between sub and supersolutions of the
Dirichlet problem when $c<0$ in $\Om$. This result is proven in
\cite{bd} under different assumptions on $F$ and $b$; thanks to
the estimate \eqref{stimau-v1}, see Remark \ref{stimadirich}, we
can show it using the same strategy of \cite{bd}, if $F$ satisfies
the conditions (F2) and (F3) and $b$ is continuous and bounded on
$\Om$.
\begin{thm}\label{dirichletcomp} Let $\Om$ be bounded. Assume that (F2) and (F3) hold, that $b$, $c$ and $g$ are continuous and bounded on
$\Om$ and $c<0$ in $\Om$. If $u\in USC(\Oms)$ and $v\in LSC(\Oms)$
are respectively sub and supersolution of
$$ F(x,Du,D^2u)+b(x)\cdot
Du|Du|^\al+c(x)|u|^{\al}u= g(x) \quad\text{in } \Om,$$ and $u\leq
v$ on $\partial\Om$ then $u\leq v$ in $\Om$.
\end{thm}
For convenience of the reader we postpone the proof of the theorem
to the next subsection.

The previous comparison result allows us to establish  the strong
minimum and maximum principles, for sub and supersolutions of the
Neumann problem even with the following more general boundary
condition $$f(x,u)+\langle Du, \overrightarrow{n}(x)\rangle=0
\quad x\in\partial\Om,$$for some
$f:\partial\Om\times\R\rightarrow\R$. We do not assume any
regularity on $f$.
\begin{prop}\label{p1}Let $\Om$ be a $C^1$ domain satisfying ($\Om$1). Assume that (F1)-(F3)
hold, that $b$ and $c$ are bounded and continuous on $\Om$ and
that $f(x,0)\leq 0$ for all $x\in\partial\Om$. If $v\in
LSC(\overline{\Om})$ is a non-negative viscosity supersolution of
\begin{equation}\label{moregencond}
\begin{cases}
F(x,Dv,D^2v)+b(x)\cdot Dv|Dv|^\al+c(x)|v|^{\al}v= 0  & \text{in} \quad\Om \\
 f(x,v)+\langle Dv, \overrightarrow{n}(x)\rangle= 0 & \text{on} \quad\partial\Om, \\
 \end{cases}
 \end{equation}
 then either $v\equiv0$ or $v>0$ in $\overline{\Om}$.
\end{prop}
\dim The assumption $(F2)$ and the fact that $F(x,p,0)=0$ imply
that
$$F(x,p,M)\geq |p|^\al \mathcal{M}_{a,A}^-M=|p| ^\al (atr(M^+)-Atr(M^-))=:H(p,M),$$ where $M=M^+-M^-$ is the minimal
decomposition of $M$ into positive and negative symmetric
matrices. It follows, since $v$ is non-negative, that it suffices
to prove the proposition when $v$ is a supersolution of the
Neumann problem for the equation
\begin{equation}\label{p1e2}
H(Dv,D^2v)+b(x)\cdot Dv|Dv|^\al-|c|_{\infty}v^{1+\al}= 0  \quad
\text{in } \Om.
 \end{equation}
Moreover we can assume $|c|_\infty>0$. Following the proof of
Theorem 2 in \cite{bd} it can be showed that $v>0$ in $\Om$. We
prove that $v$ cannot vanish on the boundary of $\Om$. We suppose
by contradiction that $x_0$ is some point in $\partial \Om$ on
which $v(x_0)=0$. For the interior sphere condition ($\Om$1) there
exist $R> 0$ and $y\in \Om$ such that the ball centered in $y$ and
of radius $R$, $B(y,R)$, is contained in $\Om$ and $x_0\in
\partial B(y,R)$. Fixed $0<\rho<R$, let us construct a subsolution
of \eqref{p1e2} in the annulus $\rho<|x-y|=r<R$. Let us consider
the function $\phi(x)=e^{-kr}-e^{-kR}$, where $k$ is a positive
constant to be determined. If we compute the derivatives of $\phi$
we get
$$D\phi(x)=-ke^{-kr}\dfrac{(x-y)}{r},\, D^2\phi(x)=\left(k^2e^{-kr}+\frac{k}{r}e^{-kr}\right)\dfrac{(x-y)\otimes(x-y)}{r^2}-\frac{k}{r}e^{-kr}I.$$
The eigenvalues of $D^2\phi(x)$ are $k^2e^{-kr}$ of multiplicity
$1$ and $-ke^{-kr}/r$ of multiplicity $N-1$. Then
\begin{equation*}\begin{split} &H(D\phi,D^2\phi)+b(x)\cdot
D\phi|D\phi|^\al-|c|_{\infty}\phi^{1+\al}\\ &\geq
e^{-(\al+1)kr}\left (ak^{\al+2}-\left
(A\dfrac{N-1}{\rho}+|b|_\infty \right )k^{\al+1}
-|c|_\infty\right).\end{split}\end{equation*} Take $k$ such that
$$ak^{\al+2}-\left (A\dfrac{N-1}{\rho}+|b|_\infty \right)k^{\al+1}
-|c|_\infty>\epsilon,$$ for some $\epsilon
>0$, then $\phi$ is a strict subsolution of the equation \eqref{p1e2}.
Now choose $m>0$ such that
$$m(e^{-k\rho}-e^{-kR})=v_1:=\text
{inf}_{|x-y|=\rho}v(x)>0,$$ and define $w(x)=m(e^{-kr}-e^{-kR})$.
By homogeneity $w$ is still a subsolution of \eqref{p1e2} in the
annulus $\rho<|x-y|<R$, moreover $w=v_1\leq v $ if $|x-y|=\rho$
and $w= 0 \leq v$ if $|x-y|=R$. Then by the comparison principle,
Theorem \ref{dirichletcomp}, $w\leq v$ in the entire annulus.

Now let $\delta$ be a positive number smaller than $R-\rho$. In
$B(x_0,\delta)\cap\overline{\Om}$ it is again $w\leq v$, in fact
where $|x-y|>R$ it is $w< 0 \leq v$; moreover $w(x_0)=v(x_0)=0$.
Then $w$ is a test function for $v$ at $x_0$. But
$$H(Dw(x_0),D^2w(x_0))+b(x_0)\cdot
Dw(x_0)|Dw(x_0)|^\al-|c|_{\infty}w^{1+\al}(x_0)> 0,$$ and
$$f(x_0,w(x_0))+ \langle Dw(x_0), \overrightarrow{n}(x_0)\rangle=f(x_0,0)+  \dfrac{\partial w}{\partial \overrightarrow{n}}(x_0)
\leq-kme^{-kR}<0.$$ This contradicts the definition of $v$.
Finally $v$ cannot be zero in $\overline{\Om}$.\finedim
\begin{rem}{\em By Proposition \ref{p1} the supersolutions in the
definition \eqref{primoaut} are positive in the whole
$\Oms$.}\end{rem}

\begin{prop}\label{p2}Let $\Om$ be a $C^1$ domain satisfying ($\Om$1). Assume that
(F1)-(F3) hold, that $b$ and $c$ are bounded and continuous on
$\Om$ and that $f(x,0)\geq 0$ for all $x\in\partial\Om$. If $u\in
USC(\overline{\Om})$ is a non-positive viscosity subsolution of
\eqref{moregencond} then either $u\equiv0$ or $u<0$ in
$\overline{\Om}$.
\end{prop}
\dim The proof is similar to the proof of Proposition \ref{p1},
observing that (F1) and the fact that $F(x,p,0)=0$ imply that
$$F(x,p,M)\leq |p|^\al(A\text{tr}(M^+)-a\text{tr}(M^-)).$$
\finedim

\begin{thm}[Maximum Principle for $c\leq 0$]\label{macpcless0}Assume the hypothesis of Proposition \ref{p2}. In addition suppose that
$\Om$ is bounded, $c\leq 0$, $c\not\equiv 0$ and $r\rightarrow
f(x,r)$ is non-decreasing on $\R$. If $u\in USC(\overline{\Om})$
is a viscosity subsolution of \eqref{moregencond} then $u\leq 0$
in $\Oms$. The same conclusion holds also if $c\equiv 0$ in the
following two cases
\begin{itemize}
\item[(i)] $\Om$ is a $C^2$ domain and there exists $\xs\in \partial\Om$ such that $\langle b(\xs), \overrightarrow{n}(\xs)\rangle>0$,
$f(\xs,r)>0$ for any $r>0$ and $S\leq 0$, where $S$ is the
symmetric operator corresponding to the second fundamental form of
$\partial \Om$ in $\xs$ oriented with the exterior normal to
$\Om$;
\item[(ii)]There exists $\xs\in \partial\Om$ such that $f(\xs,r)>0$ for any $r>0$ and $u$ is a strong subsolution.
\end{itemize}
\end{thm}
\dim Let $u$ be a subsolution of \eqref{moregencond} and
$c\not\equiv 0$. First let us suppose $u\equiv k=$const. By
definition
$$c(x)|k|^\al k\geq 0\qquad \text{in }\Om,$$ which implies $k\leq
0$.

Now we assume that $u$ is not a constant. We argue by
contradiction; suppose that $\max_{\Oms}u=u(x_0)>0$, for some
$x_0\in\Oms$. Define $\widetilde{u}(x):=u(x)-u(x_0)$. Since $c\leq
0$ and $f$ is non-decreasing, $\widetilde{u}$ is a non-positive
subsolution of \eqref{moregencond}. Then, from Proposition
\ref{p2}, either $u\equiv u(x_0)$ or $u<u(x_0)$ in
$\overline{\Om}$. In both cases we get a contradiction.

Let us turn to the case $c\equiv 0$. Suppose that $\Om$ is a $C^2$
domain, $\langle b(\overline{x}),
\overrightarrow{n}(\overline{x})\rangle>0$, $S\leq 0$ and
$f(\overline{x},r)>0$ for any $r>0$ and some point
$\overline{x}\in
\partial \Om$. We have to prove that $u$ cannot be a positive
constant. Suppose by contradiction that $u\equiv k$. In general,
if $\phi$ is a $C^2$ function, $\overline{x}\in
\partial \Om$ and $S\leq 0$ in $\xs$, then
$(D\phi(\overline{x})-\lambda \overrightarrow{n}(\xs),
D^2\phi(\overline{x}))\in J^{2,+}\phi(\overline{x})$, for $\lambda
\geq0$ (see \cite{cil} Remark 2.7). Hence $(-\lambda
\overrightarrow{n}(\overline{x}),0)\in J^{2,+}u(\overline{x})$.
But $$f(\overline{x},k)-\lam\langle
\overrightarrow{n}(\overline{x}),\overrightarrow{n}(\overline{x})\rangle=
f(\overline{x},k)-\lambda>0,$$ for $\lambda>0$ small enough, and
$$G(\overline{x},k,-\lambda \overrightarrow{n}(\overline{x}),0)=-\lambda^{\al+1}\langle b(\overline{x}), \overrightarrow{n}(\overline{x})\rangle<0.$$
This contradicts the definition of $u$.

Finally if $u$ is a strong subsolution, $f(\overline{x},r)>0$ for
 $r>0$ and some $\overline{x}\in\partial\Om$, $u\equiv k>0$,
then the boundary condition is not satisfied at $\xs$ for $p=0$.
\finedim
\begin{rem}{\em Under the same assumptions of Theorem
\ref{macpcless0}, but now with $f$ sa\-ti\-sfying $f(x,0)\leq 0$
for all $x\in\partial\Om$ and with $f(\xs,r)<0$ for any $r<0$ and
some $\xs\in\partial\Om$ in (i) and (ii), using Proposition
\ref{p1} we can prove the minimum principle, i.e., if $u\in
LSC(\Oms)$ is a viscosity supersolution of \eqref{moregencond}
then $u\geq 0$ in $\Oms$.}\end{rem}
\begin{rem}\emph{$C^2$ convex sets satisfy the condition $S\leq 0$ in
every point of the boundary.}\end{rem}
\begin{rem}{\em If $c\equiv 0$ and $f\equiv 0$ a counterexample to the maximum
principle is given by the positive constants.}\end{rem}

\subsection{The threshold for the Maximum Principle}

\noindent In this subsection and in the rest of the paper we
always assume that $\Om$ is bounded and of class $C^2$, that $F$
satisfies (F1)-(F4), that $b$ and $c$ are continuous on $\Oms$.

\begin{thm}[Maximum Principle for $\lam<\lams$]\label{maxpneum} Let $\lam <\lams$ and let $u\in
USC(\Oms)$ be a viscosity subsolution of
\begin{equation}\label{equazprincmax}
\begin{cases}
 G(x,u,Du,D^2u)+\lam |u|^{\al}u = 0  & \text{in} \quad\Om \\
 \langle Du,\overrightarrow{n}(x)\rangle= 0 & \text{on} \quad\partial\Om, \\
 \end{cases}
 \end{equation}then $u\leq 0$ in $\Oms$.
 \end{thm}
 \begin{rem}{\em Similarly it is possible to prove that if
 $\lam<\underline{\lam}$ and $v$ is a supersolution of
 \eqref{equazprincmax} then $v\geq 0$ in $\Oms$.}
 \end{rem}
\begin{cor}The quantities $\lams$ and $\underline{\lam}$ are
finite.
\end{cor}
\dim It suffices to observe that $\lams,\,\underline{\lam}\leq
|c|_{\infty}$, since when the zero order coefficient is
$c(x)+|c|_\infty$ the maximum and the minimum principles do not
hold. The theorems fail respectively for the positive and negative
constants. \finedim

In the proof of Theorem \ref{maxpneum}  the Lemma \ref{regolarita}
is one of the main ingredient. Furthermore, we need the following
two results. The first one is an adaptation of Lemma 1 of
\cite{bd} for supersolutions of the Neumann boundary value
problem; the second one is a Lemma due to Barles and Ramaswamy,
\cite{br}.
\begin{lem}\label{lemxdivy}Let $v\in LSC(\Oms)$ be a viscosity
supersolution of
\begin{equation*}
\begin{cases} F(x,Dv,D^2v)+b(x)\cdot Dv|Dv|^\al-\beta(v(x))=
g(x) & \text{in} \quad\Om \\
 \langle Dv,\overrightarrow{n}(x)\rangle= 0 & \text{on} \quad\partial\Om, \\
 \end{cases}
 \end{equation*} for some functions $g,\beta\in USC(\Oms)$. Suppose that $\xs\in \Oms$
is a strict local minimum of $v(x)+C|x-\xs|^qe^{-kd(x)}$,
$k>\frac{q}{2r}$, where $r$ is the radius in the condition
($\Om2$) and $q>\max\{2,\frac{\al+2}{\al+1}\}$. Moreover suppose
that $v$ is not locally constant around $\xs$. Then
$$-\beta(v(\xs))\leq g(\xs).$$\end{lem}

\begin{rem}\emph{Similarly, if $\beta,$ $g\in LSC(\Oms)$, $u\in USC(\Oms)$ is a
supersolution, $\overline{x}$ is a strict local maximum of
$u(x)-C|x-\xs|^qe^{-kd(x)}$, $k>\frac{q}{2r}$,
$q>\max\{2,\frac{\al+2}{\al+1}\}$ and $u$ is not locally constant
around $\xs$, it can be proved that
$$-\beta(u(\xs))\geq g(\xs).$$}\end{rem}
\begin{lem}\label{barles}If $X,Y\in S(N)$ satisfy
\begin{equation*}-\zeta\left(%
\begin{array}{cc}
  I & 0 \\
  0 & I \\
\end{array}%
\right)\leq \left(%
\begin{array}{cc}
  X & 0 \\
  0 & -Y \\
\end{array}%
\right)\leq \zeta\left(%
\begin{array}{cc}
  I & -I \\
  -I & I \\
\end{array}%
\right)
\end{equation*} then we have
 $$X-Y\leq-\frac{1}{2\zeta}(tX+(1-t)Y)^2 \quad \text{for
all }t\in[0,1].$$
\end{lem}

 \dims \textbf{of Theorem \ref{maxpneum}.}
 Let $\tau\in ]\lam,\lams[$, then by definition there exists
 $v>0$ in $\Oms$ bounded viscosity solution of
 \begin{equation}\label{soprasol1}
 \begin{cases}
 G(x,v,Dv,D^2v)+\tau v^{\al+1} \leq 0  & \text{in} \quad\Om \\
 \langle Dv,\overrightarrow{n}(x)\rangle \geq 0 & \text{on} \quad\partial\Om. \\
 \end{cases}
 \end{equation}

 We argue by contradiction that $u$ has a positive maximum in
 $\Oms$.
As in \cite{bd}, we define $\gamma':=\text{sup}_{\Oms}(u/v)>0$ and
$w=\gamma v$, with $\gamma\in (0,\gamma')$ to be determined. By
homogeneity, $w$ is still a solution of \eqref{soprasol1}. Let
$\overline{y}\in\Oms$ be such that
$u(\overline{y})/v(\overline{y})=\gamma'$. Since
$u(\overline{y})-w(\overline{y})=(\gamma'
-\gamma)v(\overline{y})>0$, the supremum of $u-w$ is strictly
positive, then by upper semicontinuity there exists
$\overline{x}\in\Oms$
 such that
$$ u(\xs)-w(\xs)=\max_{\Oms}(u-w)=m>0.$$
Clearly $u(\xs)>w(\xs)>0,$ moreover $u(\xs)\leq
\gamma'v(\xs)=\frac{\gamma'}{\gamma}w(\xs),$ from which
\begin{equation}\label{w(x)}
w(\xs)\geq \frac{\gamma}{\gamma'}u(\xs).
\end{equation}

Fix $q>\max\{2,\frac{\al+2}{\al+1}\}$ and $k>q/(2r)$, where $r$ is
the radius in the condition ($\Om2$), and define for $j\in\N$ the
functions $\phi\in C^2(\Oms\times\Oms)$ and $\psi\in
USC(\Oms\times\Oms)$ by
$$\phi(x,y)=\frac{j}{q}|x-y|^qe^{-k(d(x)+d(y))},\quad\psi(x,y)=u(x)-w(y)-\phi(x,y).$$
Let $(x_j,y_j)\in \Oms\times\Oms$ be a maximum point of $\psi$,
then $m=\psi(\xs,\xs)\leq u(x_j)-w(y_j)-\phi(x_j,y_j)$, from which
\begin{equation}\label{dis1}
\frac{j}{q}|x_j-y_j|^q\leq
(u(x_j)-w(y_j)-m)e^{k(d(x_j)+d(y_j))}\leq C,
\end{equation} where $C$ is independent
of $j$. The last relation implies that, up to subsequence, $x_j$
and $y_j$ converge to some $\zs\in\Oms$ as $j\rightarrow +\infty$.
Classical arguments show that
$$\lim_{j\rightarrow+\infty}\frac{j}{q}|x_j-y_j|^q=0,\quad\lim_{j\rightarrow+\infty}u(x_j)=u(\zs),\quad
\lim_{j\rightarrow+\infty}w(y_j)=w(\zs),$$and $$u(\zs)-w(\zs)=m.$$

\textbf{Claim 1} \emph{For $j$ large enough, there exist $x_j$ and
$y_j$ such that $(x_j,y_j)$ is a maximum point of $\psi$ and
$x_j\neq y_j$.}

Indeed if $x_j=y_j$ we have
$$\psi(x_j,x)=u(x_j)-w(x)-\frac{j}{q}|x-x_j|^qe^{-k(d(x_j)+d(x))}\leq\psi(x_j,x_j)=u(x_j)-w(x_j),$$and
$$\psi(x,x_j)=u(x)-w(x_j)-\frac{j}{q}|x-x_j|^qe^{-k(d(x)+d(x_j))}\leq\psi(x_j,x_j)=u(x_j)-w(x_j).$$ Then $x_j$ is a
minimum point for
$$W(x):=w(x)+\frac{j}{q}e^{-kd(x_j)}|x-x_j|^qe^{-kd(x)},$$and a
maximum point for
$$U(x):=u(x)-\frac{j}{q}e^{-kd(x_j)}|x-x_j|^qe^{-kd(x)}.$$
We first exclude that $x_j$ is both a strict local minimum and a
strict local maximum. Indeed in that case, if $u$ and $w$ are not
locally constant around $x_j$, by Lemma \ref{lemxdivy}
\begin{equation*}(c(x_j)+\tau)w(x_j)^{\al+1}\leq (c(x_j)+\lam )u(x_j)^{\al+1}.\end{equation*}
The same result holds if $u$ or $w$ are locally constant by
definition of sub and supersolution. The last inequality leads to
a contradiction, as we will see at the end of the proof. Hence
$x_j$ cannot be both a strict local minimum and a strict local
maximum. In the first case there exist $\delta
>0$ and $ R>\delta$ such that
\begin{equation*}\begin{split}w(x_j)&=\min_{\delta\leq|x-x_j|\leq R\atop
x\in\Oms}\left(w(x)+\frac{j}{q}|x-x_j|^qe^{-k(d(x_j)+d(x))}\right)\\&=w(y_j)+\frac{j}{q}|y_j-x_j|^qe^{-k(d(x_j)+d(y_j))},\end{split}\end{equation*}
for some $y_j\neq x_j$, so that $(x_j,y_j)$ is still a maximum
point for $\psi$. In the other case, similarly, one can replace
$x_j$ by a point $y_j\neq x_j$ such that $(y_j,x_j)$ is a maximum
for $\psi$. This concludes the Claim 1.

Now computing the derivatives of $\phi$ we get
$$D_x\phi(x,y)=j|x-y|^{q-2}e^{-k(d(x)+d(y))}(x-y)-k\frac{j}{q}|x-y|^{q}e^{-k(d(x)+d(y))}Dd(x),$$
and
$$D_y\phi(x,y)=-j|x-y|^{q-2}e^{-k(d(x)+d(y))}(x-y)-k\frac{j}{q}|x-y|^{q}e^{-k(d(x)+d(y))} Dd(y).$$
Denote $p_j:=D_x\phi(x_j,y_j)$ and $r_j:=-D_y\phi(x_j,y_j).$ Since
$x_j\neq y_j$, $p_j$ and $r_j$ are different from 0 for $j$ large
enough. Indeed
\begin{equation*}\begin{split}|p_j|,|r_j|&\geq
j|x_j-y_j|^{q-1}e^{-k(d(x_j)+d(y_j))}\left(1-\frac{k}{q}|x_j-y_j|\right)\geq
\frac{j}{2}|x_j-y_j|^{q-1}e^{-2kd_0},\end{split}\end{equation*}
where $d_0=\max_{\Oms}d(x)$. Using \eqref{sferaest}, if
$x_j\in\partial \Om$ then
\begin{equation*}
\langle p_j,\overrightarrow{n}(x_j)\rangle\geq j
|x_j-y_j|^{q}e^{-kd(y_j)}\left(-\frac{1}{2r}+\frac{k}{q}\right)>0,\end{equation*}
and if $y_j\in\partial\Om$ then
\begin{equation*}
\langle r_j,\overrightarrow{n}(y_j)\rangle\leq j
|x_j-y_j|^{q}e^{-kd(x_j)}\left(\frac{1}{2r}-\frac{k}{q}\right)<0,\end{equation*}
since $k>q/(2r)$ and $x_j\neq y_j$. In view of definition of sub
and supersolution we conclude that
$$ G(x_j,u(x_j),p_j,X)+\lam
u(x_j)^{\al+1}\geq 0\quad \text{if }(p_j,X)\in
\overline{J}^{2,+}u(x_j),$$ $$G(y_j,w(y_j),r_j,Y)+\tau
w(y_j)^{\al+1}\leq 0\quad\text{if }(r_j,Y)\in
\overline{J}^{2,-}w(y_j).$$

Applying Theorem 3.2 of \cite{cil} for any $\epsilon>0$ there
exist $X_j, Y_j\in \emph{S(N)}$ such that $(p_j,X_j)\in
\overline{J}^{2,+}u(x_j)$, $(r_j,Y_j)\in \overline{J}^{2,-}w(y_j)$
and

 \begin{equation}\label{tm1lemis}-\left(\frac{1}{\epsilon}+\|D^2\phi(x_j,y_j)\|\right)\left(%
\begin{array}{cc}
  I & 0 \\
  0 & I \\
\end{array}%
\right)\leq \left(%
\begin{array}{cc}
  X_j & 0 \\
  0 & -Y_j \\
\end{array}%
\right)\leq D^2\phi(x_j,y_j)+\epsilon (D^2\phi(x_j,y_j))^2.
\end{equation}

\textbf{Claim 2} \emph{$X_j$ and $Y_j$ satisfy
\begin{equation}\label{claimmatrix}-\zeta_j\left(%
\begin{array}{cc}
  I & 0 \\
  0 & I \\
\end{array}%
\right)\leq \left(%
\begin{array}{cc}
  X_j-\widetilde{X_j} & 0 \\
  0 & -Y_j+\widetilde{Y_j} \\
\end{array}%
\right)\leq \zeta_j\left(%
\begin{array}{cc}
  I & -I \\
  -I & I \\
\end{array}%
\right),\end{equation} where $\zeta_j=Cj|x_j-y_j|^{q-2}$, for some
positive constant $C$ independent of $j$ and  some matrices
$\widetilde{X_j},\,\widetilde{Y_j}=O(j|x_j-y_j|^{q}).$}

To prove the claim we need to estimate $D^2\phi(x_j,y_j)$.
\begin{equation*}\begin{split}D^2\phi(x_j,y_j)&=\frac{j}{q}|x_j-y_j|^qD^2(e^{-k(d(x_j)+d(y_j))})+
D(e^{-k(d(x_j)+d(y_j))})\otimes
\frac{j}{q}D(|x_j-y_j|^q)\\&+\frac{j}{q}D(|x_j-y_j|^q)\otimes
D(e^{-k(d(x_j)+d(y_j))})+e^{-k(d(x_j)+d(y_j))}\frac{j}{q}D^2(|x_j-y_j|^q).\end{split}\end{equation*}
We denote
$$A_1:=\frac{j}{q}|x_j-y_j|^qD^2(e^{-k(d(x_j)+d(y_j))}),$$
$$A_2:=De^{-k(d(x_j)+d(y_j))}\otimes
\frac{j}{q}D(|x_j-y_j|^q)+\frac{j}{q}D(|x_j-y_j|^q)\otimes
D(e^{-k(d(x_j)+d(y_j))}),$$
$$A_3:=e^{-k(d(x_j)+d(y_j))}\frac{j}{q}D^2(|x_j-y_j|^q).$$
For $A_1$ and $A_3$ we have $$A_1\leq Cj|x_j-y_j|^q \left(%
\begin{array}{cc}
  I & 0 \\
  0 & I \\
\end{array}%
\right),$$
$$ A_3\leq (q-1)j|x_j-y_j|^{q-2}\left(%
\begin{array}{cc}
  I & -I \\
  -I & I \\
\end{array}%
\right).$$ Here and henceforth, as usual, the letter $C$ denotes
various constants independent of $j$. Now we consider the quantity
$\langle A_2(\xi,\eta),(\xi,\eta)\rangle$ for $\xi,\,\eta\in\R^N$.
We have
\begin{equation*}\begin{split}\langle
A_2(\xi,\eta),(\xi,\eta)\rangle&=2kj|x_j-y_j|^{q-2}e^{-k(d(x_j)+d(y_j))}[\langle
Dd(x_j)\otimes(x_j-y_j)(\eta-\xi),\xi\rangle\\&+\langle
Dd(y_j)\otimes(x_j-y_j)(\eta-\xi),\eta\rangle]\\&\leq C
j|x_j-y_j|^{q-1}|\xi-\eta|(|\xi|+|\eta|)\\&\leq C
j|x_j-y_j|^{q-1}\left(\frac{|\xi-\eta|^2}{|x_j-y_j|}+\frac{(|\xi|+|\eta|)^2}{4}|x_j-y_j|\right)\\&\leq
C \left[j|x_j-y_j|^{q-2}|\xi-\eta|^2+
j|x_j-y_j|^{q}(|\xi|^2+|\eta|^2)\right].\end{split}\end{equation*}
The last inequality can be rewritten  equivalently in this
way$$A_2\leq C
j|x_j-y_j|^{q-2}\left(%
\begin{array}{cc}
  I & -I \\
  -I & I \\
\end{array}%
\right)+C j|x_j-y_j|^{q}\left(%
\begin{array}{cc}
  I & 0 \\
  0 & I \\
\end{array}%
\right).$$ Finally if we choose $$\epsilon
=\frac{1}{j|x_j-y_j|^{q-2}},$$ we get the same estimates for the
matrix $\epsilon(D^2\phi(x_j,y_j))^2$. In conclusion we have
\begin{equation*}\begin{split}D^2\phi(x_j,y_j)+\epsilon(D^2\phi(x_j,y_j))^2&\leq C
j|x_j-y_j|^{q-2}\left(%
\begin{array}{cc}
  I & -I \\
  -I & I \\
\end{array}%
\right)\\& +C j|x_j-y_j|^{q}\left(%
\begin{array}{cc}
  I & 0 \\
  0 & I \\
\end{array}%
\right).\end{split}\end{equation*} Hence, since
$\|D^2\phi(x_j,y_j)\|\leq Cj|x_j-y_j|^{q-2},$ \eqref{tm1lemis}
implies \eqref{claimmatrix} and the Claim 2 is proved.

\textbf{Claim 3}
\emph{$F(x_j,p_j,X_j-\widetilde{X_j})-F(y_j,r_j,Y_j-\widetilde{Y_j})\leq
o_j$, where $o_j\rightarrow0$ as $j\rightarrow+\infty$.}

First we need to know that the quantity $j|x_j-y_j|^{q-1}$ is
bounded uniformly in $j$. This is a simple consequence of Lemma
\ref{regolarita}. Indeed, since $m>0$ and $w$ is positive and
bounded, the estimate \eqref{stimau-v1} holds for $u$ and $w$;
then using it in \eqref{dis1} and dividing by $|x_j-y_j|\neq 0$ we
obtain
$$\frac{j}{q}|x_j-y_j|^{q-1}\leq C_0e^{k(d(x_j)+d(y_j))}\leq C.$$
Consequently, there exists $R>0$ such that for large $j$
\begin{equation}\label{stimapr}C\zeta_j|x_j-y_j|\leq\frac{j}{2}|x_j-y_j|^{q-1}e^{-k(d(x_j)+d(y_j))}\leq|p_j|,\,|r_j|
\leq 2j|x_j-y_j|^{q-1}\leq R.\end{equation} Denote for simplicity
$Z_j:=X_j-\widetilde{X_j}$ and $W_j:=Y_j-\widetilde{Y_j}$. By
\eqref{claimmatrix} and Lemma \ref{barles} with $t=0$, we have
$$Z_j-W_j\leq -\frac{1}{2\zeta_j}W_j^2.$$ As in the appendix of
\cite{bdl} we use the previous relation, the Cauchy-Schwarz's
ine\-qua\-li\-ty and the properties of $F$ to get the estimate of
the claim
\begin{equation*}\begin{split}F(x_j,p_j,Z_j)-F(y_j,r_j,W_j)&=F(x_j,p_j,Z_j)-F(x_j,p_j,W_j)+F(x_j,p_j,W_j)\\&-F(y_j,p_j,W_j)
+F(y_j,p_j,W_j)-F(y_j,r_j,W_j)\\&
\leq-\frac{a}{2\zeta_j}|p_j|^\al\text{tr}W_j^2+C_1|x_j-y_j|^\theta|p_j|^\al\|W_j\|\\&+C_2|p_j|^{\al-\nu}|p_j-r_j|^\nu\|W_j\|
\leq-\frac{a}{2\zeta_j}|p_j|^\al\text{tr}W_j^2\\&
+\frac{a}{4\zeta_j}|p_j|^\al\text{tr}W_j^2+\frac{C_1^2|x_j-y_j|^{2\theta}|p_j|^{2\al}\zeta_j}{a|p_j|^\al}\\&
+\frac{a}{4\zeta_j}|p_j|^\al\text{tr}W_j^2+\frac{C_2^2|p_j|^{2(\al-\nu)}|p_j-r_j|^{2\nu}\zeta_j}{a|p_j|^\al}
\\&=C\zeta_j|x_j-y_j|^{2\theta}|p_j|^{\al}+C\zeta_j|p_j|^{\al-2\nu}|p_j-r_j|^{2\nu}.
\end{split}\end{equation*}
Now consider the first term of the last quantity. Using
\eqref{stimapr} we have
$$C\zeta_j|x_j-y_j|^{2\theta}|p_j|^{\al}\leq
\frac{C\zeta_j|x_j-y_j|^{2\theta}|p_j|^{\al+1}}{\zeta_j|x_j-y_j|}\leq
C R^{\al+1}|x_j-y_j|^{2\theta-1},$$ and the last term goes to 0 as
$j\rightarrow+\infty$ since $\theta>\frac{1}{2}.$ It remains to
estimate $C\zeta_j|p_j|^{\al-2\nu}|p_j-r_j|^{2\nu}$. Observe that
$$|p_j-r_j|\leq
2k\frac{j}{q}|x_j-y_j|^{q}=C\zeta_j|x_j-y_j|^2,$$ then we have
\begin{equation*}\begin{split}C\zeta_j|p_j|^{\al-2\nu}|p_j-r_j|^{2\nu}&=
C|p_j|^{\al+1}\frac{\zeta_j}{|p_j|}\frac{|p_j-r_j|^{2\nu}}{|p_j|^{2\nu}}\leq
\frac{CR^{\al+1}}{|x_j-y_j|}|x_j-y_j|^{2\nu}\\&=CR^{\al+1}|x_j-y_j|^{2\nu-1}.\end{split}\end{equation*}
Also the last quantity goes to 0 as $j\rightarrow+\infty$ since
$\nu>\frac{1}{2}$ and this concludes the Claim 3.

Now using  the properties of $F$ and the fact that $u$ and $w$ are
respectively sub and supersolution we compute
\begin{equation*}\begin{split}
-(\lam+c(x_j))u(x_j)^{\al+1}&\leq F(x_j, p_j,X_j)+b(x_j)\cdot
p_j|p_j|^\al
\\&\leq F(x_j,p_j,X_j-\widetilde{X_j})+
b(x_j)\cdot p_j|p_j|^\al\ +|p_j|^\al O\left(j|x_j-y_j|^{q}\right)
\\&\leq F(y_j,r_j,Y_j-\widetilde{Y_j})+b(x_j)\cdot
p_j|p_j|^\al+|p_j|^\al O\left(j|x_j-y_j|^{q}\right)+ o_j\\& \leq
-(\tau+c(y_j))w(y_j)^{\al+1}+ b(x_j)\cdot p_j|p_j|^\al-b(y_j)\cdot
r_j |r_j|^\al\\&+(|p_j|^\al\vee|r_j|^\al)
O\left(j|x_j-y_j|^{q}\right)+o_j.
\end{split}\end{equation*}
 Sending $j\rightarrow+\infty$ we obtain
\begin{equation}\label{tm1ultima}-(\lam+c(\zs))u(\zs)^{\al+1}\leq
-(\tau+c(\zs))w(\zs)^{\al+1}.\end{equation} Indeed
$o_j\rightarrow0$ as $j\rightarrow+\infty$  and
$$(|p_j|^\al\vee|r_j|^\al) O\left(j|x_j-y_j|^{q}\right)\leq
C(j|x_j-y_j|^{q-1})^{\al+1}|x_j-y_j|\leq
CR^{\al+1}|x_j-y_j|\rightarrow0$$ as $j\rightarrow+\infty$.
Moreover, up to subsequence $p_j,\,r_j\rightarrow p_0\in \R^N$. If
$p_0\neq 0$ then
$$b(x_j)\cdot
p_j|p_j|^\al,\,b(y_j)\cdot r_j |r_j|^\al\rightarrow
b(\overline{z})\cdot p_0|p_0|^\al$$ and so the difference goes to
0, otherwise
$$|b(x_j)\cdot
p_j||p_j|^\al\leq |b(x_j)||p_j|^{\al+1}\rightarrow0\quad\text{as
}j\rightarrow+\infty.$$ The same result holds for $\,b(y_j)\cdot
r_j |r_j|^\al$.

If $\tau+c(\zs)>0$, from \eqref{w(x)} and \eqref{tm1ultima} we
have
$$-(\lam+c(\zs))u(\zs)^{\al+1}\leq-(\tau+c(\zs))\left(\frac{\gamma}{\gamma'}\right)^{\al+1}u(\zs)^{\al+1},$$
and taking $\gamma$ sufficiently close to $\gamma'$ in order that
$\frac{\tau\left(\frac{\gamma}{\gamma'}\right)^{\al+1}-\lam}{1-\left(\frac{\gamma}{\gamma'}\right)^{\al+1}}>
|c|_\infty,$ we get a contradiction. Finally if $\tau+c(\zs)\leq0$
we obtain
$$-(\lam+c(\zs))u(\zs)^{\al+1}\leq
-(\tau+c(\zs))w(\zs)^{\al+1}\leq
-(\tau+c(\zs))u(\zs)^{\al+1},$$once more a contradiction since
$\lam<\tau$. \finedim

\dims  \textbf{of Lemma \ref{lemxdivy}.} Without loss of
generality we can assume that $\xs=0$.

Since the minimum is strict there exists a small $\delta>0$ such
that $$v(0)<v(x)+C|x|^qe^{-kd(x)}\quad \text{for any
}x\in\Oms,\,0<|x|\leq\delta.$$ Since $v$ is not locally constant
and $q>1$ for any $n>\delta^{-1}$ there
 exists $(t_n,z_n)\in B(0,\frac{1}{n})^2\cap\Oms^2$ such that
 $$v(t_n)>v(z_n)+C|z_n-t_n|^qe^{-kd(z_n)}.$$
 Consequently, for $n>\delta^{-1}$ the minimum of the function
 $v(x)+C|x-t_n|^qe^{-kd(x)}$ in $\overline{B}(0,\delta)\cap\Oms$ is not achieved on $t_n$. Indeed
 $$\min_{|x|\leq\delta,\,x\in\Oms}(v(x)+C|x-t_n|^qe^{-kd(x)})\leq
 v(z_n)+C|z_n-t_n|^qe^{-kd(z_n)}<v(t_n).$$
Let $y_n\neq t_n$ be some point in $\overline{B}(0,\delta)\cap
\Oms$ on which the minimum is achieved. Passing to the limit as
$n$ goes to infinity, $t_n$ goes to 0 and, up to subsequence,
$y_n$ converges to some $y\in \overline{B}(0,\delta)\cap\Oms$. By
the lower semicontinuity of $v$ and the fact that 0 is a local
minimum of $v(x)+C|x|^qe^{-kd(x)}$ we have
$$v(0)\leq v(y)+C|y|^qe^{-kd(y)}\leq\liminf_{n\rightarrow +\infty}(v(y_n)+C|y_n|^qe^{-kd(y_n)}),$$
and using that $v(0)+C|t_n|^qe^{-kd(0)}\geq
v(y_n)+C|y_n-t_n|^qe^{-kd(y_n)},$ one has
$$v(0)\geq\limsup_{n\rightarrow +\infty}(v(y_n)
+C|y_n|^qe^{-kd(y_n)}).$$ Then
$$v(0)=v(y)+C|y|^qe^{-kd(y)}=\lim_{n\rightarrow +\infty}(v(y_n)+C|y_n|^qe^{-kd(y_n)}).$$
Since 0 is a strict local minimum of $v(x)+C|x|^qe^{-kd(x)}$, the
last equalities imply that $y=0$ and $v(y_n)$ goes to $v(0)$ as
$n\rightarrow +\infty$. Then for large $n$, $y_n$ is an interior
point of $B(0,\delta)$ so that the function
$$\varphi(x)=v(y_n)+C|y_n-t_n|^qe^{-kd(y_n)}-C|x-t_n|^qe^{-kd(x)}$$ is a test function for $v$ at $y_n$. Moreover, the gradient of $\varphi$
$$D\varphi(x)=-Cq|x-t_n|^{q-2}e^{-kd(x)}(x-t_n)+kC|x-t_n|^qe^{-kd(x)}Dd(x)$$
is different from 0 at $x=y_n$ for small $\delta$, indeed
$$|D\varphi(y_n)|\geq
C|y_n-t_n|^{q-1}e^{-kd(y_n)}(q-k|y_n-t_n|)\geq
C|y_n-t_n|^{q-1}e^{-kd(y_n)}(q-2k\delta)>0.$$

 Using
\eqref{sferaest}, if $y_n\in\partial \Om$ we have
\begin{equation*}\begin{split}\langle
D\varphi(y_n),\overrightarrow{n}(y_n)\rangle\leq
C|y_n-t_n|^q\left(\frac{q}{2r}-k\right)<0,
\end{split}\end{equation*} since $k>q/(2r)$. Then we conclude that
$$F(y_n,D\varphi(y_n),D^2\varphi(y_n))+b(y_n)\cdot
D\varphi(y_n)|D\varphi(y_n)|^\al-\beta(v(y_n))\leq g(y_n).$$ This
inequality together with the condition (F2) implies that
\begin{equation}\label{lem12}-|D\varphi(y_n)|^\al A\text{tr}(D^2\varphi(y_n))^-+b(y_n)\cdot
D\varphi(y_n)|D\varphi(y_n)|^\al-\beta(v(y_n))\leq
g(y_n).\end{equation} Observe that
$D^2\varphi(y_n)=|y_n-t_n|^{q-2}M,$ where $M$ is a matrix such
that $\text{tr} M^+$ and $\text{tr} M^-$ are bounded by a constant
independent of $\delta$ and $n$. Hence, from \eqref{lem12} we get
\begin{equation*}C_0|y_n-t_n|^{\al(q-1)+q-2}-\beta(v(y_n))\leq
g(y_n),\end{equation*}for some constant $C_0$, where the exponent
$\al(q-1)+q-2=q(\al+1)-(\al+2)>0$. Passing to the limit, since
$\beta$ and $g$ are upper semicontinuous we get $$-\beta
(v(0))\leq g(0),$$ which is the desired conclusion. \finedim We
conclude  sketching the proof of Theorem \ref{dirichletcomp}.

 \dims \textbf{of Theorem
\ref{dirichletcomp}.} Suppose by contradiction that
$\max_{\Oms}(u-v)=m>0$. Since $u\leq v$ on the boundary, the
supremum is achieved inside $\Om$. Let us define for $j\in\N$ and
some $q>\max\{2,\frac{\al+2}{\al+1}\}$
$$\psi(x,y)=u(x)-v(y)-\frac{j}{q}|x-y|^q.$$ Suppose that
$(x_j,y_j)$ is a maximum point for $\psi$ in $\Oms^2$. Then
$|x_j-y_j|\rightarrow0$ as $j\rightarrow+\infty$ and up to
subsequence $x_j,y_j\rightarrow\xs$, $u(x_j)\rightarrow u(\xs)$,
$v(y_j)\rightarrow v(\xs)$ and $j|x_j-y_j|^q\rightarrow0$ as
$j\rightarrow+\infty$. Moreover, $\xs$ is such that
$u(\xs)-v(\xs)=m$ and we can choose $x_j\neq y_j$. Recalling by
Remark \ref{stimadirich} that the estimate \eqref{stimau-v1} holds
in $\Om$, we can proceed as in the proof of Theorem \ref{maxpneum}
to get $$-c(\xs)|u(\xs)|^\al u(\xs)\leq -c(\xs)|v(\xs)|^\al
v(\xs).$$ This is a contradiction since $c(\xs)<0$. \finedim

\subsection{The Maximum Principle for $c(x)$ chan\-ging sign: an example.}
In the previous subsections we have proved that $G(x,u,Du,D^2u)$
with the Neumann boundary condition satisfies the maximum
principle if $c(x)\leq 0$ or without condition on the sign of
$c(x)$ provided $\lams>0$. In this subsection we want to prove
that this two cases don't coincide, i.e., that there exists some
$c(x)$ which changes sign in $\Om$ such that the associated
principal eigenvalue $\lams$ is positive. To prove this, by
definition of $\lams$, it suffices to find a function $c(x)$
changing sign for which there exists a bounded positive solution
of
\begin{equation}\label{equazsottosol}
\begin{cases}
 F(x,Dv,D^2v)+b(x)\cdot Dv|Dv|^\al+c(x)|v|^{\al}v \leq -m  & \text{in} \quad\Om \\
 \langle Dv,\overrightarrow{n}(x)\rangle\geq 0 & \text{on} \quad\partial\Om, \\
 \end{cases}
 \end{equation}where $m>0$.

In the rest of this subsection we will construct an explicit
example of such function. For simplicity, let us suppose that
$b\equiv 0$ and $\Om$ is the ball of center 0 and radius R. We
will look for $c$ such that:

\begin{equation}\label{cprop}
\begin{cases}
c(x)<0 &\text{if } R-\epsilon<|x|\leq R\\
c(x)\leq -\beta_1 & \text{if }\rho<|x|\leq R-\epsilon\\
c(x)\leq \beta_2 & \text{if }|x|\leq \rho,\\
 \end{cases}
 \end{equation}
 where $0<\rho<R-\epsilon$ and $\epsilon$, $\beta_1,\, \beta_2$ are positive constants which satisfy a
suitable inequality. Remark that in the ball of radius $\rho$
$c(x)$ may assume positive values.

In order to construct a supersolution, we define the function
\begin{equation}\label{vsubs}
v(x):=\begin{cases}
D &\text{if } R-\epsilon<|x|\leq R\\
E|x|^2-E(R+\rho-\epsilon)|x|+D+E\rho(R-\epsilon) &\text{if } \rho<|x|\leq R-\epsilon\\
D+1-e^{k(|x|-\rho)} & \text{if }|x|\leq \rho,\\
 \end{cases}
 \end{equation}
where D, E, k are positive constants to be chosen later.

\begin{lem}\label{proprietau}The function v defined in \eqref{vsubs} has the following properties
\begin{itemize}
\item[(i)] $v$ is continuous on $\overline{B}(0,R)$ and of class $C^2$ in the sets
$B(0,\rho)\setminus\{0\},\,B(0,R-\epsilon)\setminus\overline{B}(0,\rho),\,\overline{B}(0,R)\setminus
\overline{B}(0,R-\epsilon)$;
\item[(ii)]$v$ is positive provided
$D>\frac{E}{4}(R-\rho-\epsilon)^2$;
\item[(iii)] $J^{2,-}v(x)=\emptyset$ if $x=0,\,|x|=R-\epsilon$ and if $|x|=\rho$
provided $E(R-\rho-\epsilon)>k$.
\end{itemize}
\end{lem}
\dim The proof of (i) is a very simple calculation.

For (ii) we observe that $v$ is positive if $R-\epsilon\leq|x|\leq
R$ and $|x|\leq \rho$ since $D, k>0$. In the region $\{\rho\leq
|x|\leq R-\epsilon\}$ $v$ is positive on the boundary where takes
the value D, while in the interior
$Dv(x)=2Ex-E(R+\rho-\epsilon)\frac{x}{|x|}=0$ if
$|x|=\frac{R+\rho-\epsilon}{2}$. In such points
$v(x)=-\frac{E}{4}(R-\rho-\epsilon)^2+D$, then they are global
minimums where $v$ takes positive value if
$D>\frac{E}{4}(R-\rho-\epsilon)^2$.

Now we turn to (iii). Let $\widehat{x}\in\Om$ be such that
$|\widehat{x}|=\rho$ and let $(p,X)\in J^{2,-}v(\widehat{x})$,
then by definition of semi-jet
\begin{equation}\label{subsemijet}v(x)\geq v(\widehat{x})+\langle p,x-\widehat{x}\rangle
+\frac{1}{2}\langle X(x-\widehat{x}),x-\widehat{x}\rangle +
o(|x-\widehat{x}|^2),\end{equation} as $x\rightarrow \widehat{x}$.
If we take $x=\widehat{x}+t\overrightarrow{n}(\widehat{x})$, for
$t>0$, where
$\overrightarrow{n}(\widehat{x})=\frac{\widehat{x}}{|\widehat{x}|}$
is the exterior normal to the sphere of radius $\rho$ at
$\widehat{x}$, then $|x|>\rho $ and dividing \eqref{subsemijet} by
$t$ we have
$$\frac{v(\widehat{x}+t\overrightarrow{n}(\widehat{x}))-v(\widehat{x})
}{t}\geq p_n +O(t),$$ where $p_n=p\cdot
\overrightarrow{n}(\widehat{x})$. Letting $t\rightarrow 0^+$ we
get $$p_n\leq \langle
2E\widehat{x}-E(R+\rho-\epsilon)\frac{\widehat{x}}{|\widehat{x}|},\frac{\widehat{x}}{|\widehat{x}|}\rangle=-E(R-\rho-\epsilon).$$
On the other hand, if we take
$x=\widehat{x}-t\overrightarrow{n}(\widehat{x})$, $t>0$, in
\eqref{subsemijet} and divide by $-t$, letting $t\rightarrow0^+$
we get $$p_n\geq \langle
-ke^{k(|\widehat{x}|-\rho)}\frac{\widehat{x}}{|\widehat{x}|},\frac{\widehat{x}}{|\widehat{x}|}\rangle=-k.$$
In conclusion $$E(R-\rho-\epsilon)\leq -p_n\leq k.$$ Assuming the
hypothesis  in (iii) the previous condition cannot never be
satisfied, then $J^{2,-}v(\widehat{x})=\emptyset$.

In the same way it can be proved that if $\widehat{x}\in\Om$ is
such that $|\widehat{x}|=R-\rho$ and $(p,X)\in
J^{2,-}v(\widehat{x})$ then $$ E(R-\rho-\epsilon)\leq p_n\leq
0,$$and clearly also this condition cannot be satisfied,
consequently $J^{2,-}v(\widehat{x})=\emptyset$.

Finally it is easy to see that  $ J^{2,-}v(0)=\emptyset$. \finedim

\begin{prop}There exist $\epsilon,\,\beta_1,\,\beta_2>0$ such that  for any $c(x)$ satisfying \eqref{cprop} the function $v$ defined in
\eqref{vsubs} is a positive continuous viscosity solution of
\eqref{equazsottosol}.
\end{prop}
\dim Clearly $v$ satisfies the boundary condition. Since the
semi-jet $J^{2,-}v(x)$ is empty if $|x|=\rho$, $|x|=R-\epsilon$
and $x=0$, in such points we have nothing to test. In
$B(0,\rho)\setminus\{0\},\,B(0,R-\epsilon)\setminus\overline{B}(0,\rho),\,B(0,R)\setminus
\overline{B}(0,R-\epsilon)$ $v$ is of class $C^2$, then it
suffices to prove that $v$ is a classical solution of
\eqref{equazsottosol} in these open
sets. \\
 \noindent \emph{Case I: $R-\epsilon<|x|<R.$}

Since $c<0$ and continuous on $\{R-\epsilon\leq|x|\leq R\},$ we
have
\begin{equation}\label{caso1} c(x)v^{\al+1}=c(x)D^{\al+1}\leq -m_1<0.\end{equation}
Hence, by definition $v$ is supersolution.
\\
\noindent \emph{Case II: $\rho<|x|<R-\epsilon$.}

In this set
$$Dv(x)=E[2|x|-(R+\rho-\epsilon)]\frac{x}{|x|},\quad
D^2v(x)=2EI-E(R+\rho-\epsilon)\frac{1}{|x|}\left(I-\frac{x\otimes
x}{|x|^2}\right).$$Since
$-(R-\rho-\epsilon)\leq2|x|-(R+\rho-\epsilon)\leq
R-\rho-\epsilon$, using (F2) we compute
\begin{equation*}\begin{split}F&(x,Dv,D^2v)\leq
E^{\al+1}(R-\rho-\epsilon)^\al\left[\frac{2AN(R-\epsilon)-a(R+\rho-\epsilon)N+a(R+\rho-\epsilon)}{R-\epsilon}\right]\\
&=E^{\al+1}\frac{(R-\rho-\epsilon)^\al}{R-\epsilon}\{N[(A-a)(R-\epsilon)+A(R-\epsilon)-a\rho]+a(R+\rho-\epsilon)\}.\end{split}\end{equation*}
Observe that all the factors in the last member are positive.
Using the last computation, the fact that in the minimum points
$v$ takes the value $D-\frac{E}{4}(R-\rho-\epsilon)^2$ (see the
proof of Lemma \ref{proprietau})  and that $c\leq -\beta_1$, we
have
\begin{equation}\label{caso2}\begin{split}F(x,Dv,D^2v)+c(x)v^{\al+1}&\leq E^{\al+1}\frac{(R-\rho-\epsilon)^\al}{R-\epsilon}\{N[(A-a)(R-\epsilon)
+A(R-\epsilon)-a\rho]\\&+a(R+\rho-\epsilon)\}-\beta_1\left[D-\frac{E}{4}(R-\rho-\epsilon)^2\right]^{\al+1}=:-m_2.\end{split}\end{equation}
The above quantity is negative if
\begin{equation}\label{caso2'}\begin{split}D&>\frac{E}{4}(R-\rho-\epsilon)^2
+EC,\end{split}\end{equation} where
$$C:=\frac{(R-\rho-\epsilon)^{\frac{\al}{\al+1}}}{\beta_1^{\frac{1}{\al+1}}(R-\epsilon)^{\frac{1}{\al+1}}}\{N[(A-a)(R-\epsilon)
+A(R-\epsilon)-a\rho]+a(R+\rho-\epsilon)\}^{\frac{1}{\al+1}}>0.$$

\noindent\emph{Case III: $0<|x|<\rho$.}

Here we have
$$Dv(x)=-ke^{k(|x|-\rho)}\frac{x}{|x|},\quad D^2v(x)=-k^2e^{k(|x|-\rho)}\frac{x\otimes x}{|x|^2}-ke^{k(|x|-\rho)}\frac{1}{|x|}\left(I-\frac{x\otimes
x}{|x|^2}\right).$$Then
\begin{equation}\label{caso3}\begin{split}&F(x,Dv,D^2v)+c(x)v^{\al+1}\leq
-k^{\al+1}e^{(\al+1)k(|x|-\rho)}a\left(k+\frac{N-1}{|x|}\right)+\beta_2(D+1\\&-e^{k(|x|-\rho)})^{\al+1}\leq
-k^{\al+1}e^{-(\al+1)k\rho}a\left(k+\frac{N-1}{\rho}\right)+\beta_2(D+1-e^{-k\rho})^{\al+1}=:-m3.\end{split}\end{equation}
The last quantity is negative if
\begin{equation}\label{beta2}\beta_2<\frac{k^{\al+1}e^{-(\al+1)k\rho}a\left(k+\frac{N-1}{\rho}\right)}{(D+1-e^{-k\rho})^{\al+1}}.\end{equation}

Since E must satisfy the condition in (iii) of Lemma
\ref{proprietau}, we choose
\begin{equation}\label{E}E:=\frac{k}{R-\rho-\epsilon'},\end{equation} for
$\epsilon<\epsilon'<R-\rho$. Furthermore we take
\begin{equation}\label{D}D:=\frac{E}{4}(R-\rho-\epsilon)^2+EC+\epsilon=\frac{k(R-\rho-\epsilon)^2}{4(R-\rho-\epsilon')}
+\frac{kC}{R-\rho-\epsilon'}+\epsilon.\end{equation} With this
choice of D, \eqref{caso2'}  is satisfied and $v$ is positive by
(ii) of
 Lemma \ref{proprietau}. Observe that
$$D\rightarrow
k\left\{\frac{R-\rho}{4}+\left[\frac{2NAR-(N-1)a(R+\rho)}{\beta_1R(R-\rho)}\right]^{\frac{1}{\al+1}}\right\}$$
as $\epsilon,\epsilon'\rightarrow 0^+$.

Finally we can write the relation between $\beta_1$ and $\beta_2$:
\begin{equation}\label{beta1beta2}\beta_2<\frac{k^{\al+1}e^{-(\al+1)k\rho}a\left(k+\frac{N-1}{\rho}\right)}
{\left(k\left\{\frac{R-\rho}{4}+\left[\frac{2NAR-(N-1)a(R+\rho)}{\beta_1R(R-\rho)}\right]^{\frac{1}{\al+1}}\right\}+1-e^{-k\rho}\right)^{\al+1}}.
\end{equation}
Suppose that \eqref{beta1beta2} holds for some $k>0$, then we can
choose $\epsilon'>\epsilon>0$ so small that
$$\beta_2<\frac{k^{\al+1}e^{-(\al+1)k\rho}a\left(k+\frac{N-1}{\rho}\right)}{(D+1-e^{-k\rho})^{\al+1}},$$
where D is defined by \eqref{D}. Define E as in \eqref{E}, then
$v$ is a positive supersolution of \eqref{equazsottosol} with $m$
the minimum between the quantity $m_1,m_2$ and $m_3$ defined
respectively in \eqref{caso1}, \eqref{caso2} and \eqref{caso3}.
Observe that the size of $\epsilon$ is given by
\eqref{beta1beta2}.\finedim
\begin{rem} \emph{If we call $UB(\beta_2)$ the upper bound of $\beta_2$ in \eqref{beta1beta2}, we can see that if we choose $k=\frac{1}{\rho}$
 then $UB(\beta_2)$ goes to $+\infty$ as
$\rho\rightarrow 0^+$, that is, if the set where $c$ is positive
becomes small then the values of $c$ in this set can be very
large. On the contrary, for any value of $k$, if $\rho\rightarrow
R^-$ then $UB(\beta_2)$ goes to 0. Finally for any $k$ if
$\beta_1\rightarrow 0^+$, then again $UB(\beta_2)$ goes to 0. So
there is a sort of balance between $\beta_1$ and $\beta_2$. This
behavior can be explained by the following example: consider the
equation $\Delta v +c(x)v=0$ which is a subcase of our equation
and suppose that $v>0$ in $\Oms$ is a classical solution of
$\Delta v +c(x)v\leq 0$ in $\Om $, $\frac{\partial v}{\partial
\overrightarrow{n}}\geq 0$ on $\partial \Om.$ Then dividing by $v$
and integrating by part we get
\begin{equation}\label{intc}\int_{\Om}{c(x)dx}\leq-\int_{\Om}{\frac{|Dv|^2}{v^2}dx}-
\int_{\partial \Om}{\frac{1}{v}\frac{\partial v}{\partial
\overrightarrow{n}}dS}\leq 0,\end{equation} the first inequality
being strict if $\Delta v +c(x)v\not\equiv 0$. If the
supersolution is $C^2$ piecewise with $J^{2,-}v=\emptyset$ in the
non-regular points, as the one constructed before, then we can
repeat this computation in any set where $v$ is $C^2$ getting
again
$$\int_{\Om}{c(x)dx}<0.$$ }\end{rem}
\begin{rem} \emph{The construction above can be repeated for any $C^2$ domain. The
assumptions on $c$ and the supersolution $v$ can be rewritten
respectively as follows
\begin{equation*}
\begin{cases}
c(x)<0 & \text{if }d(x)< \epsilon\\
c(x)\leq -\beta_1 & \text{if }\epsilon\leq d(x)< \delta\\
c(x)\leq \beta_2 &\text{if } d(x)\geq\delta,\\
 \end{cases}
 \end{equation*}
 \begin{equation*}
v(x):=\begin{cases}
D & \text{if }d(x)<\epsilon\\
E(\delta
+\epsilon-d(x))^2+E(\delta+\epsilon)(d(x)-\epsilon-\delta)+D+E\epsilon\delta
& \text{if }\epsilon\leq d(x)< \delta\\
D+1-e^{k(\delta-d(x))} &\text{if } d(x)\geq \delta,\\
 \end{cases}
 \end{equation*}
where $0<\epsilon<\delta$ and $d(x)$ is precisely the distance
function, not one of its $C^2$ extensions. We recall some
properties of the distance function:
\begin{itemize}
    \item There exists $\mu>0$ such that d is of class $C^2$ in
    $\Om_{\mu}:=\{x\in\Oms |\, d(x)<\mu\}$ and the eigenvalues of the hessian matrix of
    $d$ at $x$ are 0 and $\frac{k_i}{1+d(x)k_i}$, $i=1,...,N-1$,
    where $k_i$ are the principal curvatures of $\partial \Om $ corresponding to the directions orthogonal to
    $\overrightarrow{n}$
     at the point $y=x-d(x)Dd(x)$;
    \item d is semi-concave in $\Om$, i.e., there exists $s_0>0$
    such that $d(x)-\frac{s_0}{2}|x|^2$ is concave;
    \item If $J^{2,-}d(x)\neq \emptyset$, $d$ is differentiable at
    $x$ and $|Dd(x)|=1$.
\end{itemize}
We choose $\delta$ so small that in $\Om_{\delta+\delta'}$ $d$ is
of class $C^2$ for some small $\delta'>0$. Then, as in previous
example, where $\delta$ was $R-\rho$, it can be proved that $v$ is
continuous on $\Oms$, positive if
$D>\frac{E}{4}(\delta-\epsilon)^2$ and of class $C^2$ on
$\Om_{\epsilon},\,\Om_{\delta}\setminus
\overline{\Om}_{\epsilon}$. Furthermore, $J^{2,-}v(x)=\emptyset$
if $d(x)=\epsilon$ and if $d(x)=\delta$ provided
$E(\delta-\epsilon)>k$.}

{\em Let $K$ be such that $|k_i(x)|\leq K $ for all $i$ and all
$x\in\partial\Om$. Then, if $\epsilon<d(x)<\delta$ we have the
following estimate
\begin{equation*}\begin{split}F(x,Dv,D^2v)+c(x)v^{\al+1}&\leq
E^{\al+1}(\delta-\epsilon)^\al\Big\{2A+[A\delta+(A-2a)\epsilon](N-1)K
\\&+[(2A-a)\delta-a\epsilon]\frac{(N-1)K}{1-\delta
K}\Big\}-\beta_1\left[D-\frac{E}{4}(\delta-\epsilon)^2\right]^{\al+1}.
\end{split}\end{equation*}
Now suppose $d(x)>\delta$, then $v(x)=D+1-e^{k(\delta-d(x))}$. Let
$\xs\in\Om$ be such that $d(\xs)>\delta$ and let $\psi$ be a $C^2$
function such that $(v-\psi)(x)\geq(v-\psi)(\xs)=0$ for all $x$ in
a small neighborhood of $\xs$. Then the function $\phi$ defined as
$$\phi(x):=-\frac{1}{k}\log(D+1-\psi(x))+\delta$$ is a $C^2$ function
in a neighborhood of $\xs$, such that
$(d-\phi)(x)\geq(d-\phi)(\xs)=0.$ This implies that
$J^{2,-}d(\xs)\neq\emptyset.$ According to some of the properties
of $d$ recalled before, on such point $d$ is differentiable,
$D\phi(\xs)=Dd(\xs)$ and $D^2\phi(\xs)\leq s_0 I$. Then it easy to
check that for $k>\frac{s_0AN}{a}$
\begin{equation*}\begin{split}F(\xs,D\psi(\xs),D^2\psi(\xs))+c(\xs)v^{\al+1}&\leq
k ^{\al+1}
e^{-(\al+1)k(R-\delta)}(s_0AN-ka)\\&+\beta_2(D+1-e^{-k(R-\delta)})^{\al+1},\end{split}\end{equation*}
where $R:=\max_{\Oms}d(x)$.}

{\em We can repeat the argument used before to conclude that $v$
is a positive solution of \eqref{equazsottosol} if $\epsilon$ is
small enough and $\beta_1$ and $\beta_2$ satisfy the following
inequality for some $k>\frac{s_0AN}{a}$
\begin{equation*}\beta_2<\frac{k^{\al+1}e^{-(\al+1)k(R-\delta)}(ka-s_0AN)}
{\left\{k\left[\frac{\delta}{4}+\left[\frac{2A+(N-1)\delta
K[A+(2A-a)(1-\delta
K)^{-1}]}{\beta_1\delta}\right]^{\frac{1}{\al+1}}\right]+1-e^{-k(R-\delta)}\right\}^{\al+1}}.
\end{equation*}}

{\em Of course the relation between $\beta_1$ and $\beta_2$ can be
bettered if we have more informations about the domain
$\Om$.}\end{rem}

\section{Some existence results} This section is devoted to the
problem of the existence of a solution of
\begin{equation}\label{existence}
\begin{cases}
 F(x,Du,D^2u)+b(x)\cdot Du|Du|^\al+(c(x)+\lam)|u|^\al u=  g(x) & \text{in} \quad\Om \\
 \langle Du, \overrightarrow{n}(x)\rangle= 0 & \text{on} \quad\partial\Om.\\
 \end{cases}
 \end{equation}
 The first existence result for \eqref{existence} is obtained when
 $\lam=0$ and $c<0$, via Perron's method. Thanks to it we will be able to prove the
 existence of a positive solution of \eqref{existence} when $g$ is
 non-positive and $\lam<\lams$, without condition on the sign of $c$. Then  we will prove the existence of a positive principal
 eigenfunction corresponding to $\lams$, that is a solution of
 \eqref{existence} when $g\equiv0$ and $\lam=\lams$. For the last
 two results we will follow the proof given in \cite{bd} for the
 analogous theorems with the Dirichlet boundary condition.

 Symmetrical results can be obtained for the eigenvalue
 $\underline{\lam}$.

 Finally, we will prove that the Neumann problem  \eqref{existence} is solvable for any right-hand side if $\lam<\min\{\lams,\underline{\lam}\}$.

 Comparison results guarantee for \eqref{existence} the uniqueness of the solution
 when $c<0$,  of the positive solution when $\lam<\lams$ and $g<0$ and of the negative solution when $\lam<\underline{\lam}$ and $g>0$.

\begin{thm}\label{princonfc<0}Suppose that $c<0$ and $g$ is continuous
on $\Oms$. If $u\in USC(\Oms)$ and $v\in LSC(\Oms)$ are
respectively  viscosity sub and supersolution of
\begin{equation}\label{princonfl<ls2}
\begin{cases}
 F(x,Du,D^2u)+b(x)\cdot Du|Du|^\al+c(x)|u|^\al u=  g(x) & \text{in} \quad\Om \\
 \langle Du, \overrightarrow{n}(x)\rangle= 0 & \text{on} \quad\partial\Om, \\
 \end{cases}
 \end{equation}
 with $u$ and $v$ bounded or $v\geq0$ and bounded, then $u\leq v$ in
$\Oms.$ Moreover \eqref{princonfl<ls2} has a unique viscosity
solution.
\end{thm}
\dim  We suppose by contradiction that $\max_{\Oms}(u-v)=m>0$.
Repeating the proof of Theorem \ref{maxpneum} taking $v$ as $w$,
we arrive to the following inequality
$$-c(\zs)|u(\zs)|^\al u(\zs)\leq -c(\zs)|v(\zs)|^\al
v(\zs),$$ where $\zs\in\Oms$ is such that $u(\zs)-v(\zs)=m>0$.
This is a contradiction since $c(\zs)<0$.

The existence of a solution follows from Perron's method of Ishii
\cite{i2} and the comparison result just proved, provided there is
a bounded subsolution and a bounded supersolution of
\eqref{princonfl<ls2}. Since $c$ is negative and continuous on
$\Oms$, there exists $c_0>0$ such that $c(x)\leq -c_0$ for every
$x\in\Oms$. Then
$$u_1:=-\left(\frac{|g|_\infty}{c_0}\right)^\frac{1}{\al+1},\quad u_2:=\left(\frac{|g|_\infty}{c_0}\right)^\frac{1}{\al+1}$$ are
respectively a bounded sub and supersolution of
\eqref{princonfl<ls2}.

Put
$$u(x):=\sup\{\varphi(x)|\,u_1\leq
\varphi\leq u_2\text{ and }\varphi \text{ is a subsolution of
\eqref{princonfl<ls2} } \},$$then $u$ is a solution of
\eqref{princonfl<ls2}. We first show that the upper semicontinuous
envelope of u defined as $$u^*(x):=\lim_{ \rho\downarrow
0}\,\sup\{u(y):y\in\Oms\text{ and } |y-x|\leq \rho\}$$ is a
subsolution of \eqref{princonfl<ls2}. Indeed if $(p,X)\in
J^{2,+}u(x_0)$ and $p\neq 0$ then by the standard arguments of the
Perron's method it can be proved that $G(x_0,u(x_0),p,X)\geq
g(x_0)$ if $x_0\in\Om$ and $(-G(x_0,u(x_0),p,X)+g(x_0))\wedge
\langle p,\overrightarrow{n}(x_0)\rangle\leq 0$ if
$x_0\in\partial\Om$.

Now suppose $u^*\equiv k$ in a neighborhood of $x_0\in\Oms$. If
$x_0\in\partial \Om$ clearly $u^*$ is subsolution in $x_0$. Assume
that $x_0$  is an interior point of $\Om$. We may choose a
sequence of subsolutions $(\varphi_n)_n$ and a sequence of points
$(x_n)_n$ in $\Om$ such that $x_n\rightarrow x_0$ and
$\varphi_n(x_n)\rightarrow k$. Suppose that $|x_n-x_0|<a_n$ with
$a_n$ decreasing to 0 as $n\rightarrow+\infty$. If, up to
subsequence, $\varphi_n$ is constant in $B(x_0,a_n)$ for any $n$,
then passing to the limit in the relation
$c(x_n)|\varphi_n(x_n)|^\al \varphi_n(x_n)\geq g(x_n)$ we get
$c(x_0)|k|^\al k\geq g(x_0)$ as desired. Otherwise, suppose that
for any $n$ $\varphi_n$ is not constant in $B(x_0,a_n)$. Repeating
the argument of Lemma \ref{lemxdivy} we find a sequence
$\{(t_n,y_n)\}_{n\in\N}\subset\Om^2$ and a small $\delta>0$ such
that $|t_n-x_0|<a_n$, $|y_n-x_0|\leq \delta$, $t_n\neq y_n$,
$\varphi_n(x)-|x-t_n|^q\leq\varphi_n(y_n)-|y_n-t_n|^q$ for any
$x\in B(x_0,\delta)$, with $q>\max\{2,\frac{\al+2}{\al+1}\}$ and
$u^*\equiv k$ in $\overline{B}(x_0,\delta)$. Up to subsequence
$y_n\rightarrow y\in \overline{B}(x_0,\delta)$ as
$n\rightarrow+\infty$. We have
\begin{equation*}\begin{split}k& =\lim_{n\rightarrow+\infty}(\varphi_n(x_n)-|x_n-t_n|^q)\leq
\liminf_{n\rightarrow+\infty}(\varphi_n(y_n)-|y_n-t_n|^q)\\&\leq\limsup_{n\rightarrow+\infty}(\varphi_n(y_n)-|y_n-t_n|^q)\leq
k-|y-x_0|^q.\end{split}\end{equation*} The last inequalities imply
that $y=x_0$ and $\varphi_n(y_n)\rightarrow k$. Then for large $n$
$y_n$ is an interior point of $B(x_0,\delta)$ and $\phi_n(x):=
\varphi_n(y_n)-|y_n-t_n|^q + |x-t_n|^q$ is a test function for
$\varphi_n$ at $y_n$. Passing to the limit as
$n\rightarrow+\infty$ in the relation $G(y_n,\varphi_n(y_n),
D\phi_n(y_n),D^2\phi_n(y_n)) \geq g(y_n)$, we get again
$c(x_0)|k|^\al k\geq g(x_0)$. In conclusion $u^*$ is a subsolution
of \eqref{princonfl<ls2}. Since $u^*\leq u_2$, it follows from the
definition of $u$ that $u=u^*$.

Finally the lower semicontinuous envelope of $u$ defined as
$$u_*(x):=\lim_{ \rho\downarrow0}\,\inf\{u(y):y\in\Oms\text{ and }
|y-x|\leq \rho\}$$is a supersolution. Indeed, if it is not, the
Perron's method provides  a viscosity subsolution of
\eqref{princonfl<ls2} greater than $u$, contradicting the
definition of $u$. If $u_*\equiv k$ in a neighborhood of
$x_0\in\Om$ and $c(x_0)|k|^\al k>g(x_0)$ then for small $\delta$
and $\rho$, the subsolution is
\begin{equation*}u_{\delta,\rho}(x):=
\begin{cases}\max\{u(x),k+\frac{\delta \rho^2}{8}-\delta|x-x_0|^2\} &\text{if }|x-x_0|<\rho,\\
u(x)&\text{otherwise}.\end{cases}
\end{equation*} Hence
$u_*$ is a supersolution of \eqref{princonfl<ls2} and then, by
comparison, $u^*=u\leq u_*$, showing that $u$ is continuous and is
a solution.

The uniqueness of the solution is an immediate consequence of the
comparison principle just proved. \finedim
\begin{thm}\label{comparisonl<ls}
Suppose $g\in LSC(\Oms)$, $h\in USC(\Oms)$, $h\leq0$, $h\leq g$
and $g(x)>0$ if $h(x)=0$. Let $u\in USC(\Oms)$ be a viscosity
subsolution of \eqref{existence} and $v\in LSC(\Oms)$ be a bounded
positive viscosity supersolution of \eqref{existence} with $g$
replaced by $h$. Then $u\leq v\quad\text{in }\Oms.$
 \begin{rem}{\em The existence of a such $v$ implies
 $\lam\leq\lams$.}\end{rem}
 \end{thm}
 \dim
 It suffices to prove the theorem for $h<g$. Indeed, for $l>1$ the function defined by $v_l:=lv$ is a
 supersolution of \eqref{existence} with right-hand side $l^{\al+1}h(x)$. By the assumptions on $h$ and $g$, $l^{\al+1}h<g$. If $u\leq lv$
for any $l>1$, passing to the limit as $l\rightarrow1^+$, one
obtains $u\leq v$ as desired.

Hence we can assume $h<g$. By upper semicontinuity
$\max_{\Oms}(h-g)=-M<0$. Suppose by contradiction that $u>v$
somewhere in $\Om$. Then there exists $\ys\in\Oms$ such that
$$\gamma':=\frac{u(\ys)}{v(\ys)}=\max_{x\in\Oms}\frac{u(x)}{v(x)}>1.$$ Define
$w=\gamma v$ for some $1\leq\gamma<\gamma'$. Since $h\leq 0$ and
$\gamma\geq1$, $\gamma^{\al+1}h\leq h$ and then $w$ is still a
supersolution of \eqref{existence} with right-hand side $h$. The
supremum of $u-w$ is strictly positive then, by upper
semicontinuity, there exists $\xs\in\Oms$ such that
$u(\xs)-w(\xs)=\max_{\Oms}(u-w)>0$. We have $u(\xs)>w(\xs)$ and
$w(\xs)\geq\frac{\gamma}{\gamma'} u(\xs)$. Repeating the proof of
Theorem \ref{maxpneum}, we get
$$g(\zs)-(\lam+c(\zs))u(\zs)^{\al+1}\leq
h(\zs)-(\lam+c(\zs))w(\zs)^{\al+1},$$where $\zs$ is some point in
$\Oms$ where the maximum of $u-w$ is attained. If $\lam+c(\zs)\leq
0$, then
$$-(\lam+c(\zs))u(\zs)^{\al+1}\leq h(\zs)-g(\zs)
-(\lam+c(\zs))w(\zs)^{\al+1}<-(\lam+c(\zs))u(\zs)^{\al+1},$$ which
is a contradiction. If $\lam+c(\zs)>0$, then
$$-(\lam+c(\zs))u(\zs)^{\al+1}\leq
h(\zs)-g(\zs)-(\lam+c(\zs))\left(\frac{\gamma}{\gamma'}\right)^{\al+1}u(\zs)^{\al+1}.$$
If we choose $\gamma$ sufficiently close to $\gamma'$ in order
that
$$|\lam+c|_\infty\left[\left(\frac{\gamma}{\gamma'}\right)^{\al+1}-1\right](\max_{\Oms}u)^{\al+1}\geq -\frac{M}{2},$$ we get
 once more a contradiction. \finedim
\begin{thm}\label{esistl<ls}
Suppose that $\lam<\lams$, $g\leq 0$, $g\not\equiv0$ and $g$ is
continuous on $\Oms$, then there exists a positive viscosity
solution of \eqref{existence}. If $g<0$, the positive solution is
unique.
\end{thm}
\dim If $\lam<-|c|_\infty$ then the existence of the solution is
guaranteed by Theorem \ref{princonfc<0}. Let us suppose $\lam\geq
-|c|_\infty$ and define by induction the sequence $(u_n)_n$ by
$u_1=0$ and $u_{n+1}$ as the solution of
\begin{equation*}
\begin{cases}
 F(x,Du_{n+1},D^2u_{n+1})+b(x)\cdot Du_{n+1}|Du_{n+1}|^\al\\ \quad+(c(x)-|c|_\infty-1)|u_{n+1}|^{\al}u_{n+1}
  =g-(\lam+|c|_\infty+1)|u_n|^{\al}u_n & \text{in} \quad\Om \\
 \langle Du_{n+1}, \overrightarrow{n}(x)\rangle= 0 & \text{on} \quad\partial\Om, \\
 \end{cases}
 \end{equation*}
which exists by Theorem \ref{princonfc<0}. By the comparison
principle, since $g\leq 0$ and $g\not\equiv 0$ the sequence is
positive and increasing.

Using the argument of Theorem 7 of \cite{bd}, thanks to the
Theorem \ref{maxpneum} and Corollary \ref{corcomp}, it can be
proved that $(u_n)_n$ is also bounded. Then, letting $n$ go to
infinity by the compactness result, the sequence converges
uniformly and, since monotone, in its whole to a function $u$
which is a solution. Moreover, the solution is positive in $\Oms$
by the strong minimum principle, Proposition \ref{p1}.

If $g<0$, the uniqueness of the positive solution follows from
Theorem \ref{comparisonl<ls}.\finedim

\begin{thm}[Existence of principal eigenfunctions]\label{esistautof}There exists $\phi>0$ in $\Oms$
viscosity solution of
\begin{equation*}
\begin{cases}
 F(x,D\phi,D^2\phi)+b(x)\cdot D\phi|D\phi|^\al+(c(x)+\lams)\phi^{\al+1}= 0 & \text{in} \quad\Om \\
 \langle D\phi, \overrightarrow{n}(x)\rangle= 0 & \text{on} \quad\partial\Om. \\
 \end{cases}
 \end{equation*}Moreover $\phi$ is Lipschitz continuous on $\Oms$.
\end{thm}
\dim Let $\lam_n$ be an increasing sequence which converges to
$\lams$. Let $u_n$ be the positive solution of \eqref{existence}
with $\lam=\lam_n$ and $g\equiv -1$. By  Theorem \ref{esistl<ls}
the sequence $(u_n)_n$ is well
 defined. Following the argument of the proof of Theorem 8 of
 \cite{bd}, we can prove that it is unbounded, otherwise one would
 contradict the definition of $\lams$.
 Then, up to subsequence $|u_n|_\infty\rightarrow+\infty$ as
 $n\rightarrow+\infty$ and defining $v_n:=\frac{u_n}{|u_n|_\infty}$
 one gets that $v_n$ satisfies \eqref{existence} with
 $\lam=\lam_n$ and $g\equiv  -\frac{1}{|u_n|_\infty^{\al+1}}$.
 Then by Corollary \ref{corcomp}, we can extract a subsequence converging to a positive
 function $\phi$ with $|\phi|_\infty=1$ which is the desired solution.
 By  Theorem \ref{regolaritathm} the solution is also Lipschitz continuous on
 $\Oms$. \finedim
 \begin{rem}\label{thmslsot}{\em With the same arguments used in the proofs of Theorems \ref{comparisonl<ls}, \ref{esistl<ls} and \ref{esistautof}
 one can prove: the comparison result between $u\in USC(\Oms)$  bounded and negative viscosity subsolution of \eqref{existence}
 and $v\in LSC(\Oms)$ supersolution of  \eqref{existence} with $g$ replaced by $h$, provided $g\geq0$, $h\leq g$ and $h(x)<0$ if $g(x)=0$;
 the existence of a negative viscosity solution of \eqref{existence}, for $\lam<\underline{\lam}$ and $g\geq 0$, $g\not\equiv0$; the
existence of a negative Lipschitz
 first eigenfunction corresponding to $\underline{\lam}$,
  i.e., a solution of \eqref{existence} with
  $\lam=\underline{\lam}$ and $g\equiv0$.}
 \end{rem}
 \begin{thm}Suppose that
$\lam<\min\{\lams,\underline{\lam}\}$ and  $g$ is continuous on
$\Oms$, then there exists a viscosity solution of
\eqref{existence}.
\end{thm}
 \dim If $h\equiv 0$, by the maximum and minimum principles the only solution is $u\equiv 0$. Let us suppose $h\not\equiv0$.
 Since $\lam<\min\{\lams,\underline{\lam}\}$ by Theorem
 \ref{esistl<ls} and Remark \ref{thmslsot} there exist
 $v_0$ positive viscosity solution of
 \eqref{existence} with right-hand side $-|g|_\infty$ and $u_0$ negative viscosity solution
 of \eqref{existence} with right-hand side $|g|_\infty$.

 Let us suppose $\lam+|c|_\infty\geq 0$. Let $(u_n)_n$ be the sequence defined in the proof of Theorem \ref{esistl<ls} with
 $u_1=u_0$, then by comparison Theorem \ref{princonfc<0} we have $u_0=u_1\leq
u_2\leq
 ...\leq v_0$. Hence, by the compactness Corollary \ref{corcomp}
 the sequence converges to a continuous function which is the
 desired solution.\finedim
\begin{ack}{\em The author wishes to thank the Professors I. Birindelli and I.
Capuzzo Dolcetta for introducing her to the problem and for
several useful discussions about the topics of this
paper.}\end{ack}

\end{document}